%
%

\magnification=1200

\font\titfont=cmr10 at 12 pt


\font\AAA=cmr14 at 12pt
\font\BBB=cmr14 at 8pt

\overfullrule=0in

\def\boxit#1{\hbox{\vrule
 \vtop{%
  \vbox{\hrule\kern 2pt %
     \hbox{\kern 2pt #1\kern 2pt}}%
   \kern 2pt \hrule }%
  \vrule}}

\def\ss{\subset}
\def\half{\hbox{${1\over 2}$}}
\def\smfrac#1#2{\hbox{${#1\over #2}$}}
\def\oa#1{\overrightarrow #1}
\def\dim{{\rm dim}}
\def\dist{{\rm dist}}

\def\log{{\rm log}}
\def\Hess{{\rm Hess}}

\def\tr{{\rm tr}}
\def\max{{\rm max}}

\def\span{{\rm span\,}}

\def\det{{\rm det}}

\def\Sym{{\rm Sym}^2}

\def\rn{\bbr^n}
\def\pp{\cp^+}
\def\plp{\cp_+}
\def\Int{{\rm Int}}
\def\cix{C^{\infty}(X)}

\def\Symn{{\Sym(\rn)}}
\def\Gpn{G(p,\rn)}
\def\fd{{\rm free-dim}}


\def\Theorem#1{\medskip\noindent {\AAA T\BBB HEOREM \rm #1.}}
\def\Prop#1{\medskip\noindent {\AAA P\BBB ROPOSITION \rm  #1.}}
\def\Cor#1{\medskip\noindent {\AAA C\BBB OROLLARY \rm #1.}}
\def\Lemma#1{\medskip\noindent {\AAA L\BBB EMMA \rm  #1.}}
\def\Remark#1{\medskip\noindent {\AAA R\BBB EMARK \rm  #1.}}
\def\Note#1{\medskip\noindent {\AAA N\BBB OTE \rm  #1.}}
\def\Def#1{\medskip\noindent {\AAA D\BBB EFINITION \rm  #1.}}

\def\Ex#1{\medskip\noindent {\AAA E\BBB XAMPLE \rm    #1.}}

\def\pf{\medskip\noindent {\bf Proof.}\ }
\def\qed{\hfill  $\vrule width5pt height5pt depth0pt$}

\def\hk{\_{\rm l}\,}
\def\n{\nabla}

   \def\cc{{\cal C}}   \def\cb{{\cal B}}  
   \def\cp{{\cal P}}

\def\cd{{\cal D}}

\def\cp{{\cal P}}
\def\cf{{\cal F}}

\def\vf{\varphi}

\def\wt{\widetilde}
\def\wh{\widehat}

\def\and{\qquad {\rm and} \qquad}

\def\bbr{{\bf R}}\def\bbh{{\bf H}}\def\bbo{{\bf O}}
\def\bbc{{\bf C}}

\def\bbz{{\bf Z}}
\def\bbp{{\bf P}}

\def\a{\alpha}
\def\b{\beta}
\def\d{\delta}
\def\e{\epsilon}
\def\f{\phi}

\def\l{\lambda}

\def\s{\sigma}
\def\x{\xi}

\def\D{\Delta}

\def\G{\Gamma}
\def\O{\Omega}

\def\psh{plurisubharmonic }

\def\lloc{L^1_{\rm loc}}

\def\bo{\partial \Omega}
\def\Ob{\overline{\O}}

\def\PSH{{ \rm PSH}}
\def\SH{{\rm SH}}

\def\MA{MA}

 \def\ppsh{$\pp$-plurisubharmonic}
\def\fp{$\pp$-plurisubharmonic }

\def\Symn{\Sym(\rn)}
 \def\ci{C^{\infty}}
\def\USC{{\rm USC}}
\def\SA{{\rm SA}}
\def\fa{{\rm\ \  for\ all\ }}
\def\ppc{$\pp$-convex}
\def\dir{Dirichlet }
\def\ft{{\wt F}}

\def\M{{\bf M}}
\def\N#1{\O_{#1}}

\def\AA{1}
\def\BB{2}

\def\EE{3}
\def\GG{4}
\def\HH{5}
\def\II{6}
\def\UU{7}
\def\DD{8}
\def\KK{9}
\def\RR{10}
\def\LL{11}
\def\JJ{12}
\def\VV{A}
\def\WW{E}

\ 
\vskip .3in

\centerline{\titfont PLURISUBHARMONICITY IN  }
\smallskip

\centerline{\titfont  A GENERAL GEOMETRIC CONTEXT}
\bigskip

\centerline{\titfont F. Reese Harvey and H. Blaine Lawson, Jr.$^*$}
\vglue .9cm
\smallbreak\footnote{}{ $ {} \sp{ *}{\rm Partially}$  supported by
the N.S.F. }

\vskip .8in
\centerline{\bf ABSTRACT} \medskip
  \font\abstractfont=cmr10 at 10 pt

{{\parindent= .6in
\narrower\abstractfont \noindent
Recently the authors have explored new concepts of plurisubharmonicity 
and pseudoconvexity, with much of the attendant analysis,
in the  context of calibrated manifolds.  Here a much broader extension is
made. This development covers a wide variety of geometric situations, including,  
for example, Lagrangian plurisubhamonicity and convexity.  It also applies
in a number of non-geometric situations.
Results include: fundamental properties of $\pp$-\psh functions, 
plurisubharmonic distributions and regularity, $\pp$-convex domains and $\pp$-convex boundaries,
topological restrictions on and construction of such domains, continuity of upper envelopes,
and solutions of the Dirichlet problem for related Monge-Amp\`ere-type equations.

}}

\vfill\eject

\centerline{\bf TABLE OF CONTENTS} \bigskip

{{\parindent= .1in\narrower\abstractfont \noindent

\qquad \AA. Introduction.\smallskip

\qquad \BB.     Geometrically Defined Plurisubharmonic Functions.    \smallskip


\qquad \EE.    More General  PSH-Functions Defined by an Elliptic cone $\pp$.    \smallskip

 \qquad \GG.    $\pp$-Plurisubharmonic Distributions.   \smallskip

 \qquad \HH.    Upper-Semi-Continuous  $\pp$-Plurisubharmonic Functions.   \smallskip

\qquad  \II.   Some Classical Facts that Extend to $\pp$-Plurisubharmonic Functions.  \smallskip

\qquad   \UU.  The Dirichlet Problem -- Uniqueness.  \smallskip

\qquad   \DD.  The Dirichlet Problem -- Existence.  \smallskip

\qquad   \KK.   $\pp$-Convex Domains.   \smallskip

\qquad   \RR.   Topological Restrictions on $\pp$-Convex Domains.   \smallskip

\qquad    \LL.   $\plp$-Free Submanifolds.  \smallskip

\qquad   \JJ.    $\pp$-Convex Boundaries.  \smallskip

\vskip .3in

 
 \quad Appendix \VV:  The Maximum Principle and Subaffine Functions.\smallskip
   
 \quad  Appendix B:  Hessians of Plurisubharmonic Distributions. \smallskip
 
 \quad  Appendix C:   Convex Elliptic Sets in $\Symn$.\smallskip
 
 \quad  Appendix D:   The Dirichlet Problem for  Convex Elliptic Sets.\smallskip
 
 \quad  Appendix E:    Elliptic $\MA$-operators / G\aa rding-Hyperbolic Polynomials on Sym$^2(\rn)$.\smallskip

}}

\vfill\eject

\centerline{\bf \AA. Introduction.}
\medskip

  Recently the authors have shown that the concepts of plurisubharmonicity 
and pseudoconvexity from complex analysis carry over, along with many of the basic results,
to other geometries,  including calibrated and symplectic geometry.
In this paper the same ideas and results are extended to a broad geometric context.
The core concept is that of an {\sl elliptic cone}.  This is a closed convex cone $\pp$
in the space $\Symn$ of symmetric $n\times n$-matrices, with the property that
the relative interior of its polar dual $\plp$ consists of positive definite matrices.

A function $u$ of class $C^2$ on an open set $X \ss \rn$ is   defined  to be 
{\sl $\pp$-plurisubharmonic} if $\Hess_x u \in \pp$ at every point $x$.

Basic geometric examples are constructed as follows.  Fix an integer $p$,  $1\leq p \leq n$, and
denote by $G(p,\rn)$ the Grassmannian of $p$-planes in $\rn$. Embed
$$
G(p,\rn)\ \ss\ \Symn
$$
by associating to each $p$-plane $\x$, the orthogonal projection $P_\x :\rn \to \x\ss\rn$.
Now let $G\ss G(p,\rn)$ be any compact subset, and define $\plp(G)$ (note the {\sl lower} plus) to be the closed convex cone
in $\Symn$ generated by  $G$.  Then a function $u\in C^2(X)$ is $\pp(G)$-plurisubharmonic if and only if
$$
\tr_\x \{\Hess_x u\} \ \geq\ 0\qquad\forall\, x\in X \ \ {\rm and}\ \ \forall \, \x\in G
$$
where $\tr_\x A \equiv \langle A, P_\x\rangle$ denotes the trace of $A$ on the $p$-plane $\x$.

Important examples of this type are where $G=G(\phi)$ consists of the $p$-planes associated
to a calibration $\phi$ of degree $p$ (such as the K\"ahler, or Special Lagrangian, or Associative,
Coassociative or Cayley calibrations).  Other interesting cases are where $G$ is the set of 
all Lagrangian $n$-planes in $\bbc^n$, or where $G=G(p,\rn)$.

This geometric case has the following interesting feature.  A function $u\in C^2(X)$ is 
$\pp(G)$-\psh if and only if its restriction to every minimal $G$-submanifold of $X$ is subharmonic
in the induced metric. (A $G$-{\sl submanifold} is a $p$-dimensional submanifold of $X$ all
of whose tangent planes lie in $G$.)

Of course the concept of an elliptic cone is much broader than the geometric case.  Nevertheless,
a surprising bulk of classical pluripotential theory carries over to this context.  The notion of 
$\pp$-plurisubharmonicity extends from $C^2$-functions to distributions, and every such 
distribution is actually in $L^1_{\rm loc}$ and has a unique upper semi-continuous 
representative with values in
$[-\infty, \infty)$.  The set $\PSH(X)$ of such   functions has all the classical properties. 
For example, if $u,v\in \PSH(X)$, then $\max\{u,v\}\in\PSH(X)$. Also, $\PSH(X)$ is closed under decreasing limits and uniform limits.  An important fact is that if $\cf\ss \PSH(X)$ is a family which is locally 
bounded above, then (the upper semicontinuous regularization) of $\sup_{v\in \cf} v$ is in $\PSH(X)$.
This enables one to apply the Perron process.   

There is a notion of $\pp$-convexity generalizing the concept of pseudo-convexity in complex 
analysis.  Given a compact set $K\ss X$, we define its $\pp$-convex hull to be the set $\wh K$
of points $x$ with 
$$
u(x)\ \leq\ \sup_K u \qquad \ \ {\rm for\ all \ smooth\ } u\in \PSH(X).
$$
Then $X$ is said to be {\sl $\pp$-convex} if for all $K\ss\ss X$ we have $\wh K\ss\ss X$.
 It is proved that {\sl  $X$ is $\pp$-convex if and only if $X$ admits a strictly $\pp$-plurisubharmonic
exhaustion function.}

Given a compact domain $\O\ss X$ with smooth boundary $\bo$, there is also a notion of $\pp$-convexity (and strict $\pp$-convexity) of the boundary.  It is shown that if $\bo$ is strictly
$\pp$-convex, then $\O$ itself is $\pp$-convex.

There is also a concept which generalizes the notion from complex geometry of being totally real.
In \S 10 we introduce  the notion of a linear subspace $V\ss\rn$ which is {$\plp$-free}.
In the geometric case this means that $V$ contains no $G$-planes, that is, there
are no $\x\in G$ with $\x\ss V$. Then the {\sl free dimension of $\plp$}, denoted fd$(\plp)$,  is  defined to be the
largest dimension of a $\plp$-free subspace of $\rn$, and we have the following generalization
of the   Andreotti-Frankel Theorem.

\Theorem{\RR.5} {\sl Any $\pp$-convex domain has the homotopy type of a CW-complex
 of dimension $\leq $ \ {\rm fd}$(\plp)$.}
\medskip

The integer fd$(\plp)$ is often easily computable, particularly in the geometric cases.
See \S \RR\ for examples.

 A submanifold is said to be {\sl $\pp$-free} if all of its
tangent planes are $\plp$-free. This extends the notion of totally 
real submanifolds in complex geometry.   In geometric cases any submanifold
of dimension $\leq p$ is free. Generic submanifolds of dimension $\leq $ fd$(\plp)$ are
  $\plp$-free on an open dense subset.  Therefore,   examples of $\plp$-free submanifolds are
  easy to construct.  This leads to lots of $\pp$-convex domains via
   the following analogue of the Grauert Tubular Neighborhood
  Theorem.

\Theorem{\LL.4}  {\sl Suppose $M$ is a $\plp$-free closed submanifold of $X\ss \rn$.  
Then there exists a fundamental neighborhood system $\cf(M)$ of $M$ 
consisting of $\pp$-convex domains.  Moreover,
\smallskip
\item{a)}  $M$ is a deformation retract of each $U\in \cf(M)$.

\smallskip
\item{b)} Each compact subset $K\ss M$ satisfies $K={\wh K}_U$ for  all  $U\in \cf(M)$.
}
\medskip

The methods used in [HW$_{1,2}$]  to generalize the Grauert Theorem
 extend  to prove this very general   result.

Freeness of submanifolds and convexity of their tubular neighborhoods are 
related by the following fact.  Let $M$ be a closed submanifold of an open
subset $X\ss \rn$.  Then {\sl $M$ is $\plp$-free if and only if the square of the distance
to $M$ is strictly $\pp$\psh at each point of $M$} (and hence in a neighborhood of $M$).
More generally we have the following result.

\Theorem{\LL.3}  {\sl Consider the two classes of closed sets.
\smallskip

\item{1)}  Closed subsets $Z\ss M$ of a $\plp$-free submanifold $M\ss X$.

\smallskip

\item{2)}  Zero sets $Z=\{f=0\}$ of non-negative strictly  $\pp$-\psh functions $f$. \smallskip

\noindent
Locally these two classes are the same.}
 
 \medskip
 
 One of the main  results of this paper is the existence and uniqueness of solutions
 to the Dirichlet Problem for functions which are {\sl  $\pp$-taut} or {\sl $\pp$ partially
 pluriharmonic}.   For functions which are $C^2$ this means that $\Hess_x u\in \partial
 \pp$ for all $x\in X$.  More generally for $u\in \PSH(X)$ this notion is defined via a duality
 involving the  subaffine functions, which are discussed in Appendix A. The main results are the following.

\Theorem{\DD.1. (The Dirichlet Problem -- Existence)}  {\sl
Suppose  $\O$ is a bounded domain in $\rn$ with a  strictly  $\pp$-convex boundary.
Given $\vf \in C(\bo)$, the function $u$ on $\overline\O$ defined by taking the upper envelope:
$$
\qquad u(x)\ =\ \sup \{v(x) : v\in \pp(\vf)\} \qquad\qquad{\rm where}
$$
$$
\pp(\vf)\ \equiv  \bigl\{v \ : \ v \in \USC(\Ob),\ \  v\bigr|_{\O} \in  \PSH(\O)
\ \ {\rm and\ \ }v\bigr|_{\bo}\leq \vf    \bigr\}
\eqno{(\DD.1)}$$
satisfies:
\medskip

1)\ \ \ $u\in C(\Ob)$,

\medskip

2)\ \ \ $u$ is  $\pp$ partially pluriharmonic on $\O$, 

\medskip

3)\ \ \ $u\bigr|_{\bo} = \vf$  on  $\bo$.}

\medskip

\Theorem{\UU.1.  (The Dirichlet Problem--Uniqueness)} {\sl  Suppose
$\pp$ is an elliptic cone and that $K$ is a compact subset of $\rn$.
If $u_1, u_2 \in C(K)$ are $\pp$-partially  pluriharmonic   on $\Int K$, then }
$$
u_1\ =\ u_2\ \  {\rm on}\  \partial K \qquad \Rightarrow \qquad
u_1\ =\ u_2\ \  {\rm on}\  K 
$$

\medskip

Many of the results in this paper have been subsequently generalized by the authors.
For example, in [HL$_4$] Theorems   \RR.5, \UU.1 and  \DD.1 above have been established
for fully non-linear,  degenerate elliptic equations which are purely of second order.
This paper makes extensive use of subaffine functions
and a certain duality intrinsic to these second order problems.
Subaffine functions are introduced here in Appendix A .  They play an important
role in the proof of the Uniqueness Theorem \UU.1 above.
This paper also treats the Dirichlet Problem for {\sl all branches} of the real, complex and
quaternionic Monge-Ampere equations.

In [HL$_5$] results are extended to closed subsets $F\ss J^2$ of the 2-jet bundle of functions
on $\rn$.  Here $F$ depends on all the classical variables $(x,r,p, A)\in X\times \bbr\times\rn\times
\Symn$. A Notion  of $F$-subharmonic functions is given and   all the  good properties
discussed above are established. Many of the theorems here are carried over.

In [HL$_6$]  the parallel discussion is carried out on riemannian manifolds where there
are may interesting geometric applications.
While these latter papers largely subsume the results here, we feel that this article
has valuable features. The exposition is less technical.  The cases covered here 
include many of  basic geometric interest.  Finally, since the basic sets $\pp$ are convex cones,
we are able to use convolutions and classical distribution theory.  This makes the analytic
part of the paper more widely accessible.  
The latter papers use other more technical analytic methods.  The article [HL$_4$] employs
deep results of Slodkowski [S]. The  papers [HL$_{5,6}$] employ the powerful Viscosity 
approach pioneered by Crandall, Ishii, Lions, Evans, Jensen and others (cf. [CIL], [C]).

 \vfill\eject

 \noindent
 {\bf Conventions:}  \smallskip
 
1.  Throughout this paper $X$ shall denote a connected open subset of $\rn$. We note that 
almost all of the analysis done here carries over to much more general riemannian manifolds $X$.
 \smallskip
 
 2.  Whenever $C\subset V$ is a convex cone in a finite dimensional vector space $V$ we
 shall denote by $\Int C$ the interior of $C$ in the vector subspace $W=\span C$.
 
 \vfill\eject


\centerline{\bf \BB. Geometrically Defined Plurisubharmonic Functions }
\medskip

In this section we  discuss a notion of plurisubharmonicity, for $\ci$-functions, based on a distinguished subset $G$ of the Grassmannian. We shall begin with some definitions and notation.
Let $G(p,\rn)$ denote the Grassmannian of unoriented $p$-planes through the origin in $\rn$.
Let $\Symn$ denote the vector space of quadratic forms (functions) on $\rn$.  We identify 
$G(p,\rn)$ with a subset of $\Symn$ by associating to each $\x\in G(p,\rn)$ the quadratic form
$P_\x$ corresponding to orthogonal projection of $\rn$ onto $\x$. 
The natural inner product on  $\Sym(\bbr^n)$ is given by the trace:
$
\langle A, B\rangle  \ =\ \tr AB.
$
Let $\cp$ denote the set of non-negative quadratic forms, $A\geq 0$.
This  is a closed convex cone with vertex at the origin in $\Symn$.
The interior,  ${\rm Int } \cp$, consists of the positive definite quadratic forms, $A>0$.
The extreme rays in $\cp$ are generated by the rank-1 projections $G(1,\rn)$.

The {\sl polar} of a closed convex cone $\cc$ with vertex  at the origin is defined by
$$
 \cc^0\ \equiv\ {\rm polar}\,\cc\ \equiv\ \{A : \langle A,B\rangle \geq 0\ \ {\rm for\ all\ \ }B\in \cc\}.
\eqno{(\BB.1)}
$$
The Bipolar Theorem states that $(\cc^0)^0=\cc$. Note that the cone $\cp$ is self-polar, that is 
$\cp^0=\cp$, since $A\geq0$ if and only if $\langle A,P_\x\rangle \geq0$ for all $\x\in G(1,\rn)$.
(If $x\in\rn$ is a unit vector  and $\x$  is the line through $x$, 
then $\langle A,P_\x\rangle =\langle A x,x\rangle$.)

Given $\x\in G(p,\rn)$ and $A\in \Symn$, the {\sl $\x$-trace of $A$}, defined by
$$
\tr_\x A\ =\ \langle A,P_\x\rangle \ =\ \tr\left( A\bigr|_\x   \right),
\eqno{(\BB.2)}
$$
is central to our development.


Given a  function $u\in C^\infty(X)$,  its hessian    
at a point  $x\in X$  will be denoted by $\Hess_x \, u$. This  is a quadratic form  on $\rn$, i.e., $\Hess_x\, u\in\Symn$.

\Def{\BB.1}  Suppose $G$ is a non-empty closed subset of $G(p,\rn)$.  A  function $u\in\cix$  is called  {\sl  $G$-\psh} if 
$$
\tr_\x\left(\Hess_x \, u\right)\ \geq \ 0\qquad{\rm for \  each\ } \x\in G, x\in X.
\eqno{(\BB.3)}
$$
Let $\PSH^\infty(X,G)$ denote this space of $G$-\psh functions.\medskip

Suppose $W$ is an affine $p$-plane through $x$ with tangent space $TW=\x$.  Then 
$$
\tr_\x\left(\Hess_x \, u\right)\ =\ \tr\left(\Hess_x \, u \bigr|_\x\right)\ =\ \tr\left(\Hess_x \, u \bigr|_W\right)
\ =\ \D\left(\ u \bigr|_W\right).
\eqno{(\BB.4)}
$$
 Call $W$ an {\sl affine $G$-plane} if $TW=\x \in G$. Then (\BB.3) is equivalent to the following.
$$
\D\left(\ u \bigr|_{X\cap W}\right)\ \geq\ 0\qquad {\rm for\ each\  affine\ }G{\rm  -plane\ } W .
\eqno{(\BB.3)'}
$$
That is, the restriction of $u$ to each affine $G$-plane $W$ is subharmonic.

A submanifold $M$ of $\rn$ is a {\sl $G$-submanifold} if $T_xM\in G$ at each point $x\in M$.

\Theorem{\BB.2} {\sl Suppose $M\subset X$ is a $G$-submanifold which is minimal. For each 
$u\in C^\infty(X)$ which is $G$-plurisubharmonic, the restriction of $u$ to $M$ is subharmonic in the induced 
riemannian metric on $M$.}

\pf
Recall the classical fact (cf. \S 1 in  [HL$_2$])  that if $u\in C^\infty(X)$, then 
for a minimal submanifold $M$,  the Laplace Beltrami operator of $M$ is given at $x\in M$ by 
$$\qquad\quad
\D_M\left(u\bigr|_M\right) \ =\ \tr\left\{\Hess_x \,  u\bigr|_{T_xM}\right\} \ 
=\ \tr_{T_xM} \left\{\ \Hess_x \, u \bigr|_{T_xM}  \right\}. \qquad\quad   \vrule width5pt height5pt depth0pt 
$$

\medskip
\noindent
{\bf Partially  Pluriharmonic Functions.}  In tandem with the concept of $G$-plurisubharmonicity
it is natural to define a function $u\in \cix$ to be {\sl $G$-pluriharmonic} if  
$$
\tr_\x \Hess_x u \ =\ 0  \qquad {\rm for \  each \ } \x\in G {\ \rm and \  each\ }x\in X.
\eqno{(\BB.5)}
$$
That is, $u$ is  $G$-pluriharmonic if and only if the restriction of $u$ to each affine $G$-plane is harmonic.  As in the proof of Theorem \BB.2, if $M$ is a $G$-submanifold which is minimal
and $u$ is $G$-pluriharmonic,
then $u\bigr|_M$ is harmonic in the induced riemannian metric on $M$.  Unfortunately, with rare exceptions,
the space of $G$-pluriharmonic functions is small (finite dimensional).  See the examples below.

A weakening of the definition of $G$-pluriharmonicity provides a much larger class.

\Def{\BB.3}   Suppose $u\in \PSH^\infty (X, G)$.  Then 
\smallskip
\item{1)}    $u$ is called {\sl partially $G$-pluriharmonic} if for each $x\in X$, the trace
$$
\tr_\x \Hess_x u \ =\ 0\qquad{\ \rm for\ some\ }\x\in G.
$$ 
\item{2)}   $u$ is called {\sl strictly $G$-plurisubharmonic} if for each $x\in X$, 
$$
\tr_\x \Hess_x u \ >\ 0\qquad{\ \rm  for\ all\ }\x\in G.
$$


\noindent{\bf Examples.}
There are many  geometrically interesting cases of $G$-\psh functions
to which our general theory will apply. This wealth of examples is one of the important
features of  this paper.  

A rich source   is the theory of calibrations [HL$_{1,2}$].  
Let $\phi$ be a constant coefficient $p$-form on $\rn$ with the property that
$\phi(\x)\leq 1$ for all $\x\in G(p,\rn)$.  Then we define the $\f$-Grassmannian
to be the set 
$$
G(\f)\ =\ \{\x\in G(p,\rn) : \phi(\x)=1\}
$$
In the following examples all but numbers 1,3, 13 can be constructed this way.

\medskip

\noindent
1.\  $G=G(1,\rn)$.\   $\PSH(X,G)$ is the set of {convex functions} on $X$.

\medskip

\noindent
2.\  $G=G(n,\rn) = \{I\}$ with $I\in \Symn$  the identity. \ $\PSH(X,G)$ is the set of {subharmonic  functions} on $X$.

\medskip

\noindent
3.\  $G=G(p,\rn)$ for $1<p<n$.   $\PSH(X,G)$  is called the  set of {\sl real $p$-\psh functions} on $X$.
The defining  property is  that they are subharmonic on every  affine $p$-plane.

\medskip

\noindent
4.\  $G= \bbp^{n-1}(\bbc) = G_\bbc(1,\bbc^n) \subset G(2, \bbr^{2n})$   gives the set of standard \psh functions in complex analysis.
\medskip

\noindent
5.\  $G= \bbp^{n-1}(\bbh) = G_\bbh(1,\bbh^n) \subset G(4, \bbr^{4n})$ gives the set of quaternionic  \psh functions on quaternionic $n$-space $\bbh^n$ (cf.  [Al], [AV]).
\medskip

\noindent
6.\  $G = G_\bbc(p,\bbc^n)$ for $1<p<n$  gives complex $p$-\psh functions on  $\bbc^n$.
\medskip

\noindent
7.\  $G = G_\bbh(p,\bbh^n)$ for $1<p<n$  gives  quaternionic $p$-\psh functions on  $\bbh^n$.
\medskip

\noindent
8.\  $G = \{x_1$-axis$\} \subset G(1,\bbr^n)$  gives  the {\sl horizontally convex functions}, i.e., the functions which are convex in the variable $x_1$.
\medskip

\noindent
9.\  $G = {\rm SLAG} \subset G(n,\bbc^n)$,  the set of (unoriented) special Lagrangian 
 $n$-planes in $\bbc^n$.
 \medskip

\noindent
10.\  $G = {\rm ASSOC} \subset G(3,\bbr^7)$,  the set of (unoriented) associative 3-planes in 
Im$\bbo \cong \bbr^7$, the imaginary octonions.
 \medskip

\noindent
11.\  $G = {\rm COASSOC} \subset G(4,\bbr^7)$,  the set of (unoriented) coassociative 4-planes in 
Im$\bbo$.
 \medskip

\noindent
12.\  $G = {\rm CAY} \subset G(4,\bbr^8)$,  the set of (unoriented) Cayley 4-planes in 
the octonions  $\bbo\cong\bbr^8$.
 \medskip

\noindent
13.\  $G = {\rm LAG} \subset G(n,\bbc^n)$,  the set of   Lagrangian 
 $n$-planes in $\bbc^n$.
 \medskip

\Remark{\BB.4} As noted in the introduction, for expository reasons the discussion in
this paper is confined to $\bbr^n$ with $G$ parallel. However, all of the examples above
can be carried over to general riemannian manifolds equipped with some additional structure.
Note for example that 4,6 and 13 make sense on any symplectic manifold with a compatible riemannian metric.   A  quite general  analysis  on riemannian manifolds is carried out in 
[HL$_6$].

\vskip .3in


\centerline{\bf  Elliptic Subsets $G$ of the Grassmannian.}
\medskip

In this section a notion of ellipticity is discussed which puts a very natural restriction on the subsets
$G\subset \Gpn$.

Let $\plp(G)$ denote the closed convex cone in $\Symn$, with vertex at the origin, determined by the compact set $G \subset \Symn$.  Let $\pp(G)$ denote the polar of $\plp(G)$.
Note that since $\cp=\plp(G(1,\rn))$ contains all the Grassmannians $\Gpn$, we have
$$
\plp(G)\ss \cp, \qquad{\rm and\ hence\ }\qquad  \cp\ss \pp(G),
$$
 for any $G\ss \Gpn$.  Set 
$$
S(G) =\span G = \span \plp(G).
$$
As one can see from the examples, $S(G)$ is usually a proper vector subspace of $\Symn$, 
and, in particular, $\plp(G)$ has no interior in $\Symn$.  However, considered as a subset of 
$S(G)$, the interior of $\plp(G)$ has  closure equal to $\plp(G)$. By $\Int \plp(G)$ we shall always mean
the interior of $\plp(G)$ in $S(G)$, not in $\Symn$.  In particular, $\Int\plp(G)$ is never empty.

\Def{\BB.5}  A closed subset $G\subset G(p,\rn)$ is   {\sl elliptic} 
if each $A\in\Int\plp(G)$ is positive definite.
\medskip

The following conditions on a closed subset $G\ss \Gpn$ are equivalent.
\smallskip

1)\ \ Given $x\in\rn$, if $x\hk\x=0$ for all $\x\in G$, then $x=0$.\smallskip

2)\ \ For each unit vector $e\in \rn$, $P_e$ is never orthogonal to $S(G)=\span G$.\smallskip

3) \ \ There does not exist a hyperplane $W\ss\rn$ with $G\ss \Sym(W)$.
\medskip

\noindent
To see that 1) and 2) are equivalent, note that $\langle P_e,P_\x\rangle = |e\hk \x|^2$.  If $e\perp W$, then $G\ss\Sym(W)$ if and only if $e\hk\x=0$ for all $\x\in G$, so that 2) $\Leftrightarrow$ 3).

\Def{\BB.6}  A closed subset $G\subset G(p,\rn)$ is said {\sl to involve all the variables in $\bbr^n$} 
if one of the equivalent conditions 1), 2), 3) holds.

 \Prop{\BB.7}  {\sl Suppose $G$ is a closed subset of $\Gpn$.  Then $G$ is elliptic if and only if 
 $G$ involves all of the variables in $\rn$.}
 
 \pf  If $G$ does not involve all the variables in $\rn$, then, by 3), $G$ and so also  $\plp(G)$
 are contained in $\Sym(W)$ for some hyperplane $W$.  
 This excludes the possibility that there exists an $A\in S(G)$ which is positive definite.

If $G$ involves all the variables in $\rn$, then, by 2), we have the following.  Under the orthogonal decomposition
$$
 P_e\ =\ E_e+S_e\qquad{\rm with\ } S_e\in S(G)\ \ {\rm and\ \  } E_e\perp S(G),
\eqno{(\BB.6)}
$$
the component $S_e$ is never zero.  Now choose $A\in\Int \plp(G)$.  
Since $S_e\in S(G)$, it follows that   for 
small $\e>0$ we have $A-\e S_e\in \Int \plp(G)\ss \cp$.
 Therefore, $0\leq  \langle P_e,A -\e S_e \rangle =  \langle P_e,A \rangle -\e|S_e|^2$ proving 
 that $ \langle P_e,A \rangle>0$ for all unit $e\in \rn$, i.e. $A>0$. \qed
 
 \medskip
 
 Each $A\in\Symn$ determines a constant coefficient linear second-order operator
 $\langle \Hess\, u,A \rangle$, which is elliptic if and only if $A>0$ (positive definite).
 If $A>0$,  then
 $$\D_Au =  \langle \Hess\, u,A \rangle$$
 will be called the {\sl $A$-Laplacian}.

 \Def{\BB.8} 
 Suppose $G$ is elliptic.  Then for each $A\in \Int \plp(G)$, the $A$-Laplacian
 $\D_A$ will be called a {\sl mollifying Laplacian} for $G$-\psh functions.

\medskip
\noindent
{ {\AAA M\BBB OLLIFYING \AAA L\BBB EMMA \rm  \BB.9.}}
 {\sl Suppose $G$ is elliptic and $u\in \cix$.  Then $u$ is $G$-\psh 
 if and only if $u$ is $\D_A$-subharmonic for each mollifying Laplacian $\D_A$.}
 
 \pf
 This follows from the fact that $G\ss \plp(G)$ and that $\plp(G)$ is the closure of $\Int \plp(G)$.\qed

\vfill\eject


\centerline{\bf \EE. More General Plurisubharmonic Functions Defined by an Elliptic Cone $\pp$}
\medskip

The basic properties of geometrically defined \psh functions remain valid in much greater generality.
Suppose $\pp$ is a closed convex cone in $\Symn$ with vertex at the origin.  Let $\plp$ denote the 
polar cone.  Let $S(\plp)$ denote the span of $\plp$, and let $\Int \plp$ denote the relative interior
of $\plp$ in the vector subspace $S(\plp)$ of $\Symn$. 

\Def{\EE.1}  
\smallskip
\item{1)} \ $\pp$ is a {\sl positive cone} if each $A\in \plp$ is positive, i.e. $A\geq0$.
 
\smallskip
\item{2)} \ $\pp$ is an {\sl elliptic cone} if each $A\in \Int\plp$ is positive definite.

\medskip\noindent
{\AAA R\BBB EMARK.}  Of course in the geometric case $\plp=\plp(G)$, the {\sl positivity condition}
$\plp\ss \cp$ is automatic.
\medskip

If $\pp$ is an elliptic cone (and, to  a lesser extent, if $\pp$ is a positive cone), it is reasonable to investigate $\pp$-\psh functions, 
even though they have no direct geometric interpretation.

\Def{\EE.2}  A function $u\in \cix$ is {\sl $\pp$-\psh} if   
$$
\Hess_x u \in \pp\qquad{\rm for\ each\  }x\in X.
\eqno{(\EE.1)}
$$
{\AAA R\BBB EMARK.}
If $\Hess_x u\in \partial \pp$, then $u$ is {\sl partially $\pp$-pluriharmonic}.
 If $\Hess_x u\in \Int \pp$, then $u$ is {\sl strictly $\pp$-plurisubharmonic}.
 Finally, if $\Hess_x u \perp S(\plp)$, then $u$ is {\sl  $\pp$-pluriharmonic}.

\medskip
A main point is that the Mollifying Lemma remains valid.
 
\medskip
\noindent
{\AAA M\BBB OLLIFYING \AAA L\BBB EMMA \rm  \EE.3.}
{\sl Suppose $\pp$ is an elliptic cone and $u\in \cix$.  Then $u$ is $\pp$-\psh 
if and only if $u$ is $\D_A$-subharmonic for each mollifying Laplacian.}

\Remark{\EE.4}  There is an analogue of (\BB.3)$'$.  Let $G$ denote the extreme points in the compact convex base $\cb_+ = \plp\cap\{\tr=1\}$. Then $u$ is $\pp$-\psh at $x$  if and only if 
$$
\langle  \Hess_x\, u, A\rangle\ =\ \tr\left\{(\Hess_x\,u)A\right\} \ \geq0\qquad\forall A\in G.
$$ 
However, this is not particularly interesting or useful unless the extreme points of the base
 $\cb_+=\plp\cap\{\tr=1\}$ are known.  It is easy to see in the geometric case where $G\ss\Gpn$ and $\plp=\plp(G)$, that the set of extreme points of $\plp\cap\{\tr=1\}$ is exactly $G$.
 
 \bigskip


\centerline{\bf  Reformulating Ellipticity for a Cone $\pp$.}
\medskip

The positivity and  ellipticity conditions have useful reformulations.   
First, the

\medskip
\noindent
{\bf Positivity Condition:} \ \ $\plp\ss\cp$, that is, every $A\in \plp$ is $\geq0$
\medskip
 
  can be stated in the equivalent dual form:
\medskip
\noindent
{\bf Positivity Condition:} \ \ $\cp\ss\pp$, that is, every $A\geq0$ belongs to $\pp$.
\medskip

In terms of functions $u$, this says that each convex function is $\pp$-\psh.

As noted, it is unusual for $S(\plp)=\span \plp$ to be all of $\Symn$.  However, there is a
different kind of incompleteness that should be ruled out. Suppose $e$ is a unit vector in $\rn$ and $W$ is the orthogonal
hyperplane in $\rn$.  Then $\Sym(W)$ can be considered a subspace of $\Symn$.  We say that 
$\pp$ {\sl can be defined using the variables in $W$} if 
$$
\pp\ =\ \Sym(W)^\perp\oplus \left(\pp\cap \Sym(W)\right).
\eqno{(\EE.2)}
$$
We say that $\plp$ {\sl only involves the variables in $W$} if 
 $$
\plp\ \ss\ \Sym(W).
\eqno{(\EE.2)'}
$$
 It is easy to see that (\EE.2) and (\EE.2)$'$ are equivalent.

\medskip
\noindent
{\bf Completeness Condition:} The cone $\pp$ can not be defined 
using the variables in a proper subspace $W\subset  \rn$, or equivalently, 
$\plp$ involves all the variables in $\rn$.
 
 \Prop{\EE.5}  {The cone $\pp$ is elliptic if and only if the positivity condition and the completeness condition are both satisfied.}
 
 \pf
 First note that if $\pp$ is elliptic, then $\plp=\overline{\Int\plp}\ss\cp$, i.e., the Positivity Condition is satisfied.  The Completeness Condition must also be satisfied, since $\plp\ss\Sym(W)$ excludes the possibility of $\plp$ containing $A>0$.
 
 The following fact is basic to our discussion.
 
 If the Positivity Condition   $\plp\ss\cp$ is satisfied, then for each $A\in \plp$ and $W=e^\perp$:
 $$
\langle A, P_e \rangle \  =\  A(e, e)=0, \ \  {\rm if\  and\  only\ if \ \ }  A\in  \Sym(W)\ss\Symn
 \eqno{(\EE.3)}
 $$
This follows because if $A\geq 0$ and $A(e,e)=0$, then
 $0\leq A(te+u, te+u) = 2tA(e,u)+A(u,u)$ for all $t\in \bbr$ and all $u\in W= e^\perp$. Hence,
 $A(e,u)=0$ for all $u\in W$, i.e.,  $ A\in\Sym(W)$.
 
  If  $\plp$ involves all the variables in $\rn$and the Positivity Condition $\plp\ss\cp$ is satisfied,  
  then because of (\EE.3),
 $\langle A, P_e \rangle $ cannot vanish for all $A\in \plp$,
  i.e., $S_e$ is never zero. (Recall the decomposition (\BB.6).)  This forces $A\in \Int\plp$ to be positive definite exactly as in the last paragraph of the proof of Proposition \BB.7.\qed
 \medskip
 
 \noindent
 {\bf Smoothing Maxima.}  
 As we shall discuss, many of the facts concerning classical subharmonic functions on
 $\rn$ extend, once we have a suitable definition of (non-smooth) $\pp$-\psh functions.
 However, limiting the discussion to smooth $\pp$-\psh functions, there are still several
 interesting facts that extend. One of the most basic is smoothing the maximum of two 
 $\pp$-\psh functions.  Let $M(t) \equiv \max\{t_1,...,t_m\}$ for  $t\in\bbr^m$.
 Suppose $\vf\in C^{\infty}_{\rm cpt}(\bbr^m)$, $\vf\geq0$, $\int\vf = 1$, with $\vf(-t)=\vf(t)$ and 
 spt$\,\vf\ss \{t:|t|\leq1\}$.
 
 Since $M$ is a convex function, $M_\e = M*\vf_\e $ is convex and decreases to $M$
 as $\e\to 0$.  Also, $\sum_{j=1}^m {\partial M\over \partial t_j} =1$ implies 
 $\sum_{j=1}^m {\partial M_\e\over \partial t_j} =1$, or equivalently, $M(t+\l e) = M(t)+\l$ implies
 $M_\e(t+\l e) = M_\e(t)+\l$, where $e=(1,...,1)$.  Moreover, $M_\e(t)-\e\leq M(t)\leq M_\e(t)$. 
 Finally note that ${\partial M \over \partial t_j}\geq 0$ implies that 
 ${\partial M_\e\over \partial t_j}\geq 0$.
 
 Given $u^1,...,u^m\in \cix$  and $\Psi$ a smooth function of $m$ variables, the chain rule implies
 that
 $$
 \Hess \Psi(u^1,...,u^m) \ =\ \sum_{j=1}^m {\partial \Psi  \over \partial t_j} \Hess u^j 
 +  \sum_{i,j=1}^m {\partial^2 \Psi  \over \partial t_i\partial t_j} \nabla u^i\circ \nabla u^j
 \eqno{(\EE.4)}
 $$
 \smallskip
 
 \noindent
 {\bf   Maxima Property:}  Suppose $\pp$ is a positive cone (not necessarily elliptic). Given
 $u^1,...,u^m\in \PSH^\infty(X)$, one has that:
 \smallskip
 \item{1)}  \ \ $M_\e(u^1,...,u^m) \in \PSH^\infty(X)$,
 
 \smallskip
 \item{2)}  \ \ $M_\e(u^1,...,u^m) - \e  \leq M(u^1,...,u^m) \leq M_\e(u^1,...,u^m)$,
 
 \smallskip
 \item{3)}  \ \  $M_\e(u^1,...,u^m)$ decreases to $M(u^1,...,u^m)$   as $\e\to 0$.    \medskip
 
 \pf
 Properties 2) and 3) are properties of $M(t)$ and $M_\e(t)$.  
 To prove 1) consider a more general function $\Psi$ and apply (\EE.3).
 The value of the quadratic form 
 $B\equiv \sum {\partial^2 \Psi  \over \partial t_i\partial t_j} \nabla u^i\circ \nabla u^j$
 on $\x=(\x_1,...,\x_n)\in\rn$ is $\sum {\partial^2 \Psi  \over \partial t_i\partial t_j} 
 \langle \nabla u^i, \x\rangle   \langle \nabla u^j, \x\rangle$, which is $\geq 0$ if $\Psi$ is convex.
 If each ${\partial \Psi  \over \partial t_j}\geq0$ and $\sum_j {\partial \Psi  \over \partial t_j}=1$,
 then the quadratic form $A = \sum_j {\partial \Psi  \over \partial t_j} \Hess u^j$ is a convex combination of the quadratic forms $\Hess u^j$, $j=1,...,m$.
 These assumptions are valid for $\Psi=M_\e$.
 The convexity of $\pp$ and the positivity condition imply that 
 $\Hess \Psi(u^1,...,u^m) = A+B\in \pp$, which proves 1).\qed

 \medskip
 \noindent
 {\AAA E\BBB XERCISE } \EE.6.  Suppose $\psi'\geq 0$ and $\psi''\geq0$. Show that $u\in
 \PSH^\infty(X)$ implies $\psi(u)\in\PSH^\infty(X)$.


\vskip .3in

\centerline{\bf \GG.  $\pp$-Plurisubharmonic Distributions. }
\medskip

Throughout this section we assume that  $\pp\subset \Symn$ is an elliptic cone 
(with vertex at the origin).

\Def{\GG.1}  A distribution $u\in\cd'(X)$ is  {\sl \fp}  if 
$$
(\Hess\, u )(\varphi A)\geq0 \quad {\rm  for\  all\ \ } A\in\plp
\eqno{(\GG.1)}$$
and all test functions 
$\varphi\in C^\infty_{\rm cpt}(X)$ with $\varphi\geq0$,

\medskip

\noindent
It is easy to see  Definition \GG.1 is compatible with Definition   \BB.1 for $u\in \cix\ss \cd'(X)$.

 \Note{}  Let $\PSH^{\rm dist}(X)$ denote the set of $u\in \cd'(X)$ which are \fp distributions on $X$. Obviously $\PSH^{\rm dist}(X)$ is a closed convex cone in $\cd'(X)$.
\medskip

The condition (\GG.1) for distributional $\pp$-plurisubharmonicity can be modified as follows.
The test function $\vf$ can be eliminated since we have
$$
(\Hess\,  u)(\vf A) \ =\ (\langle \Hess\,  u, A\rangle)(\vf).
$$
where $\Hess\,  u$ is a symmetric matrix with entries in $\cd'(X)$.  Therefore,
for a given $A\in \Symn$, condition (\GG.1) is equivalent to the statement that
$$
\D_A u\ =\ \langle \Hess\,  u, A\rangle\ \geq\ 0 \ \ \ {\rm is\  a\ non-negative \ measure.}
\eqno{(\GG.1)'}
$$

The Mollifying Lemma \EE.3  extends to distributions.

\medskip
\noindent
{\AAA M\BBB OLLIFYING \AAA L\BBB EMMA \rm  \GG.2.}   {\sl  Suppose $\pp$ is an elliptic cone.  A distribution $u\in\cd'(X)$ is \fp if and only if  $u$ is a $\D_A$-subharmonic distribution for each mollifying Laplacian $\D_A$ (i.e., each $A\in \Int \pp$).}

\pf
This is essentially the equivalence of (\GG.1) and  (\GG.1)$'$.  Also note that each $A\in\plp$
can be approximated by elements in $\Int\plp$.\qed
\vskip .3in

\noindent{\bf $G$-Plurisubharmonic Distributions}. 
Assume that $\pp=\pp(G)$ is geometrically defined by an elliptic subset 
$G$ of the grassmannian $G(p,\rn)$.  For each $\x\in G$, consider the degenerate Laplacian
defined by $A=P_\x$, i.e.,
$$
\D_\x u \ =\ \langle \Hess\,  u, P_\x\rangle
\eqno{(\GG.2)}
$$
Equation (\BB.3)$'$ has an extension  from $u\in \cix$ to $u\in \cd'(X)$.
\Prop{\GG.3}  {\sl Suppose $u\in \cd'(X)$.  Then $u\in \PSH(X)$ if and only if }
$$
\D_\x u      \     \geq\ 0\qquad {\sl \  for\ all \ }\x\in G.
\eqno{(\GG.3)}
$$
\pf 
Each $A\in \Int\plp(G)$  is a finite  positive linear combination of projections
$P_\x$ with $\x\in G$. Hence, (\GG.3) implies that $\D_A u \geq 0$ for each $A\in 
\Int\plp(G)$,  so that $u\in \PSH^{\rm dist}(X)$ by the mollifying Lemma \GG.2.
Conversely, if $u\in \PSH^{\rm dist}(X)$, then $\D_A u \geq 0$ for each $A\in \Int\plp(G)$.
 If $\x\in G$ and $t>0$, then for $A'\in \Int\plp(G)$, one has $A= P_\x+tA'\in \Int\plp(G)$
 since ${1\over t} P_\x+A' \in \Int\plp(G)$ for $t$ large.
 Hence, $\D_\x u + t\D_{A'} u \geq 0$ for all $t>0$, which proves that $\D_\x u\geq 0$ if $\x\in G$.
 \qed

\medskip

Many of the classical facts about subharmonic distributions immediately carry over to 
\fp distributions because of the Mollifying Lemma \GG.2.
We  list these classical facts in  \S \II.
\vskip.3in



\centerline{\bf \HH.  Upper-Semi-Continuous $\pp$-Plurisubharmonic Functions. }
\medskip

Throughout this section we assume that $\pp\ss\Symn$ is an elliptic cone (with vertex
at the origin).  Let $\USC(X)$ denote the space of $[-\infty,\infty)$-valued function on $X$
which are upper-semi-continuous, and not $\equiv -\infty$ on any component of $X$.

\Def {\HH.1}  A function $u\in \lloc(X)$ is called {\sl $\lloc$-upper-semi-continuous} if the essential limit
superior 
$$
\wt{u} (x) \ =\  {\rm ess}     \limsup_{y\to x} \, u(y)
\eqno{(\HH.1)}
$$
satisfies the conditions:
\smallskip
(i)  \ \ $\wt u \in \USC(X)$, and

(ii) \ $\wt u$ lies in the   $\lloc(X)$-equivalence class of $u$.
\medskip

\Prop{\HH.2}  {\sl  Each $\pp$-\psh distribution $u$ is $\lloc$-upper-semi-continuous.}
\medskip\noindent
The proof of the proposition will be given below.

\Def{\HH.3}  If $u\in \PSH^{\rm dist}(X)$, the associated canonical representative 
$\wt u\in \USC(X)$ is said to be an {\sl upper-semi-continuous $\pp$-\psh function.}
Let $\PSH^{\rm u.s.c.}(X)$ denote the space of 
upper-semi-continuous $\pp$-\psh functions on $X$.

\Cor{\HH.4}  {\sl The map sending $u\in \PSH^{\rm dist}(X)$ to $\wt u\in \PSH^{\rm u.s.c.}(X)$
is  an isomorphism.}
\pf
The map is surjective by definition, and injectivity is obvious.\qed

\medskip

We shall denote these equivalent spaces $\PSH^{\rm dist}(X) \cong  \PSH^{\rm u.s.c.}(X)$
simply by $\PSH(X)$ when no confusion will arise.

Classical potential theory applies to each Laplacian $\D_A$ with $A$ positive definite.  Since
$\D_A$ is obtained from the standard Laplacian $\D$ by a linear change of coordinates,
any result for the standard Laplacian $\D$ that is independent of choice of linear coordinates
applies to each $\D_A$ as well.

Let $\SH_A^{\rm u.s.c.}(X)$ denote the space of classical $\D_A$-subharmonic functions.
That is, $u\in \SH_A^{\rm u.s.c.}(X)$ if $u\in \USC(X)$ and for each compact subset $K$ of
$X$ and each $\D_A$-harmonic function $h$ on a neighborhood of $K$,
$$
u\leq h\ \ {\rm on\ \ } \partial K \quad {\rm implies}\quad u\leq h \ \ {\rm on\ \ } K
\eqno{(\HH.2)}
$$

Let $\SH_A^{\rm dist}(X)$ denote the space of   $\D_A$-distributions on $X$.
That is, $u\in \SH_A^{\rm dist}(X)$ if $u\in \cd'(X)$ and $\D_A u\geq0$ is a non-negative 
regular Borel measure on $X$.

For the standard Laplacian $\D$ on $\rn$ there are many references for the fact that 
$\SH^{\rm dist}(X)$ and $\SH^{\rm u.s.c.}(X)$ can be identified. 
More specifically, with $A=I$:
\smallskip

1)\ \ $u\in \SH_A^{\rm dist}(X)$ implies  $u\in \lloc(X)$.
\smallskip

2)\ \ $u\in \SH_A^{\rm dist}(X)$ implies  that $\wt u\in\SH_A^{\rm u.s.c.}(X)$
and $\wt u$ lies  in the $\lloc(X)$-class of $u$.
\smallskip

3)\ \ $u\in \SH_A^{\rm u.s.c.}(X)$ implies  $u\in \lloc(X)$.
\smallskip

4)\ \ $u\in \SH_A^{\rm u.s.c.}(X)$ implies  $\D_Au\geq 0$.
\smallskip

\noindent
Note that 1) and 2) provide an injective map $\SH^{\rm dist}(X) \to  \SH^{\rm u.s.c.}(X)$ 
given by $u\mapsto \wt u$, while 3) and 4) assert the surjectivity of this map.

These properties 1)---4) carry over to any $A>0$ by the appropriate linear coordinate change
on $\bbr^n$.  This proves that 
$$
\SH_A^{\rm dist}(X) \ \cong\ \SH_A^{\rm u.s.c.}(X).
\eqno{(\HH.3)}
$$

The $\lloc$-upper-semi-continuity Condition 2) can be proved as follows for $\D_A=\D$.
Let $B_r(x)$ denote the ball of radius $r$ about $x$ and $|B_r(x)|$ the volume of $B_r(x)$.
By the mean value inequality
$$
\wt u(x)\ \leq \ {1\over |B_r(x)|} \int_{B_r(x)} u\ \leq \ {\rm ess} \sup_{B_r(x)} u\ \leq \ 
\sup_{B_r(x)} \wt u.
\eqno{(\HH.4)}
$$
Since $\wt u$ is u.s.c., we have $\limsup_{y\to x} \wt u(y) = \wt u (x)$ forcing the essential lim sup to equal $\wt u(x)$.\qed

\medskip

Stated differently,  we have shown that if $u\in\cd'(X)$ is both $\D_{A_1}$-subharmonic and 
$\D_{A_2}$-subharmonic, then the two classical representatives $\wt{ u_1}, \wt{ u_2} \in \USC(X)$ are
equal.  Thus there is no ambiguity in the u.s.c. function $\wt u$ representing a $\pp$-\psh function.

\medskip
\noindent
{\bf Proof of Proposition \HH.2.}  If $u\in \PSH^{\rm dist}(X)$, then for some $A>0$, $u\in 
 \SH^{\rm dist}_A(X)$ and Condition 2) is valid. \qed
 \medskip
As a corollary, the Mollifying Lemma can be stated for u.s.c. $\pp$-\psh functions.

\medskip
\noindent
{ {\AAA M\BBB OLLIFYING \AAA L\BBB EMMA \rm  \HH.5.}} 
{\sl A function $u\in \USC(X)$ is u.s.c. $\pp$-\psh if and only if $u$ is
 u.s.c. $\D_A$-subharmonic for each mollifying Laplacian $\D_A$, i.e., each $A\in\Int \plp$.}

\medskip


\vskip .2in
\noindent {\bf Upper-Semi-Continuous $G$-Plurisubharmonic Functions.}
Suppose that $\pp=\pp(G)$ is geometrically defined by an elliptic subset $G$ of the grassmannian
$G(p,\rn)$.  Theorem \BB.2 about $C^\infty$ $G$-\psh functions, has only a weak extension 
(Proposition \GG.3)  to $G$-\psh distributions.  However, it has a strong extension to upper-semi-continuous $G$-\psh functions.

 \Theorem{\HH.6}  {\sl  Suppose $\pp$ is geometrically defined by an 
 elliptic subset $G$ of $G(p,\rn)$.
Let $u\in \PSH^{\rm u.s.c.}(X)$  
and  suppose $W$ is an affine $G$-plane with $W\cap X$ connected.
Then either   $u\bigr|_{W\cap X} \equiv -\infty$ or
\smallskip
\centerline{  $u\bigr|_{W\cap X}$ is subharmonic. }\smallskip

\noindent
 More generally,  suppose $M$ is any connected   $G$-submanifold of $X$, i.e.,  
 $T_xM\in G$ for all $x\in M$.  If  $M$ is  a minimal submanifold, then 
either $u\bigr|_{M}\equiv-\infty$ or\smallskip

\centerline{ 
$u\bigr|_{M}$ is subharmonic. }\smallskip

\noindent
in the induced riemannian metric on $M$.}

\pf
Assume that $u\in \PSH^{\rm u.s.c.}(X)$ and $u$ is not $\equiv -\infty$ on $M$.  As noted in \S \II, there exists a sequence $\{u_j\}$ of smooth $\pp$-\psh 
  functions on $X$ decreasing to $u$.  By Theorem \BB.2 each $u_j\bigr|_M$ is subharmonic. 
  Hence, the 
  decreasing limit   $u\bigr|_M$ is subharmonic.\qed

Theorem \HH.6 has a converse.

\Prop{\HH.7}  {\sl  Suppose that $u$ is a $[-\infty,\infty)$-valued u.s.c. function on a ball $B\subset\rn$
with the property that
for every affine $G$-plane $W$ in $\rn$,
either $u\bigr|_{W\cap B}\equiv-\infty$ or $u\bigr|_{W\cap B}$ is subharmonic.  If $u\in\lloc(B)$, then $u\in   \PSH^{\rm u.s.c.}(B)$.}

\pf
It suffices to show that $u\in \PSH^{\rm dist}(B)$ by  Corollary \HH.4.  By Proposition \GG.3 it suffices to show that $\D_\x u\geq0$ for each $\x\in G$.
Choose coordinates so that $\x$ is the first axis $p$-plane in $\rn$ and $(x,y)$ belongs 
to $\bbr^p\times \bbr^{n-p}=\rn$.  It suffices to show that $\int_{\rn} u \D_x \vf \geq 0$ for all
$\vf \in C^{\infty}_{\rm cpt}(\rn)$, $\vf\geq0$. Now $U(y)\equiv \int_{\bbr^p} u(x,y) \D_x \vf(x,y)
\in \lloc(\bbr^{n-p})$, and $U\geq 0$ a.e. by hypothesis. Hence,   by Fubini's Theorem
$\int_{\rn}  u\D_x \vf =\int_{\bbr^{n-p}} U(y) dy
\geq0$.


\vskip .3in

\centerline{\bf  \II.  Some Classical Facts that Extend to $\pp$-Plurisubharmonic Functions.}\smallskip

In this  section we list other useful properties of $\PSH(X)$-functions.

Some of the  standard results for $\D_A$-subharmonic functions immediately extend to 
$\pp$-\psh functions by the Mollifying Lemma \HH.5. Other facts require more discussion. In what follows, $u\in\PSH(X)$ is always the canonical, u.s.c. representative.

\medskip
\noindent
{\bf Facts that follow immediately from the Mollifying Lemma.}
\medskip

\item {(1)}  (Maxima)  \ \ $\max\{u^1,...,u^N\} \in \PSH(X)$ if $u^1,...,u^N\in
 \PSH(X)$ .
 \medskip

\item {(2)}  If $\psi$ is a convex non-decreasing function, then $\psi(u)  \in
 \PSH(X)$ for each $u\in \PSH(X)$
  \medskip

\item {(3)}   (Maximum Principle)  If $K$ is a compact subset of $X$ and $u\in \PSH (X)$,
then 
$$
u(x)\ \leq \ \sup_{\partial K} u\qquad {\rm for\ all\ \ }x\in K.
$$

\item {(4)}   (Decreasing Limits) If $\{u_j\}_{j=0}^\infty$ is a decreasing (i.e., $u_j\geq u_{j+1}$)  sequence of functions  in $\PSH (X)$ and $X$ is connected, then unless 
$u=\lim_{j\to\infty} u_j$
is identically $-\infty$, one has $u\in \PSH (X)$ and $\{u_j\}$ converges to $u$ in $\lloc(X)$.
  
 \medskip

\item {(5)}   (Increasing Limits) Suppose  $\{u_j\}_{j=0}^\infty$ is an increasing (i.e., $u_j\leq u_{j+1}$)  sequence of functions  in $\PSH (X)$. 
If the limit $u=\lim_{j\to\infty} u_j$ is locally bounded above, then the u.s.c. regularization
$u^*(x) = \limsup_{y\to x} u(y)$ of $u$ belongs to $\PSH (X)$ with 
$u^*=u$ a.e. and $\{u_j\}$ converging to $u$ in $\lloc(X)$.

 \medskip

\item {(5)$'$}   (Families Locally Bounded Above) Suppose $\cf\ss \PSH (X)$ is 
locally uniformly bounded above.  Then the upper envelope $v=\sup_{f\in\cf} f$ has u.s.c. regularization
$v^*\in \PSH (X)$ and $v^*=v$ a.e..  Moreover, there exists a sequence $\{u_j\}\ss\cf$
with $v^j = \max\{u_1,...,u_j\}$ converging to $v^*$ in $\lloc(X)$.

  
 \medskip

\item {(6)}   (Viscosity Plurisubharmonic) $u\in \PSH(X)$ if and only if $u\in \USC(X)$ and for each point $x\in X$ and each function $\vf \in C^2$ near $x$ with $u - \vf$ having a local maximum  at $x$, one has
$$\Hess_x \vf \in \pp.$$

\medskip
\noindent
{\bf Facts that do not  follow immediately from the Mollifying Lemma.}
\medskip

\item {(7)}    For each  $u\in  \PSH^{\rm u.s.c.}(X)$, there exists a decreasing sequence
of smooth functions $\{u_j\}\in  \PSH^{\infty}(X_j)$  with  $u=\lim_{j\to\infty}u_j$, where $X_j
=\{x\in X: \dist(x, \partial X) \geq 1/j \}$.

\medskip

\item {(8)} If $u^1,...,u^m \in \PSH(X)$ have the property that 
$\Hess u^j -\Lambda$ is $\pp$-positive, where $\Lambda : X\to \Symn$ is continuous,
then $\Hess\{M_\e(u^1,...,u^m)\}-\Lambda$ is $\pp$-positive.

\medskip

\item {(9)} (Richberg)  Suppose $u\in C(X)\cap \PSH(X)$ has the property that  $\Hess u-\Lambda$ 
is $\pp$-positive on $X$ where  $\Lambda : X\to \Symn$ is continuous. 
 Given $\l\in C(\overline X)$, $\l>0$ on $X$, there exists $\wt u \in \cix\cap\PSH(X)$ with 
$$
u\ \leq \  \wt u\ \leq\ u+\l \qquad {\rm on}\ \ \O
$$
such that  $\Hess \wt u - (1-\l)\Lambda$ is $\pp$-positive on $X$.

\bigskip
\noindent
{\bf Some Comments:}
\smallskip

By the ``classical case of (k)"  we will mean statement (k) with $\PSH(X)$ replaced by 
$\SH_A(X)$ with $A >0$.  
\smallskip
\noindent  {\bf (5)$'$:}  The classical case of (5)$'$ follows  from the classical case of (1) and (5) because of Choquet's Lemma, which says that for any family $\cf\ss\USC(X)$ which is uniformly bounded above, there exists a sequence $\{u_j\}\ss\cf$ such that the upper envelopes 
$v(x) = \sup_{f\in \cf} f(x)$ and $u(x) = \sup_j u_j(x)$
satisfy  the inequalities $u\leq v\leq v^*\leq u^*$ which forces $u^*=v^*$..
Note that (5)$'$ also follows directly from (1) and (5) by using Choquet's Lemma.

\smallskip
\noindent  {\bf (6):}  There is an $\e$-strict version of (6).  See Definition \JJ.6. 
See a) and b) below.

\smallskip
\noindent  {\bf (7):}   This statement can be proved as follows.
 If $u_\e =\vf_\e * u$ is a convolution smoothing, 
then $\D_A u_\e=\vf_\e *( \D_A u)$ so that each  $u_\e$ is in $\PSH^{\infty}$ on a subset of $X$ a distance $\e$ away from the boundary. 
If $I\in \plp$, then  the
convolutions $u_\e =\vf_\e * u$ with $\vf_\e(x) = \e^{-n}\vf({|x|\over \e})$ based on a radial function
$\vf(|x|)$, decrease monotonically  to $u$ as $\e\to 0$.  Since $\D_A$ is equivalent to $\D$ under a linear coordinate change, we can also find $\vf$ such that  $u_\e = \vf_\e * u\searrow u$  if
 $u$ is $\D_A$-subharmonic.
\qed

\smallskip
\noindent {\bf (8) and (9):}  A matrix of distributions, such as $\Hess u-\Lambda$, is defined to be
 {\sl $\pp$-positive} if $\langle \Hess u-\Lambda, A\rangle\geq 0$ is a non-negative measure
 for all $A\in\plp$. The proofs of (8) and (9) are the same as  in the several complex variable case.
   See Richberg [R] and [D] 
 Lemma 5.18 e) for (8) and Theorem 5.21 for (9).

\medskip
\noindent
{\bf  Pluriharmonicity and Strict Plurisubharmonicity.}
It is straightforward to extend the definition   of pluriharmonicity to distributions.

\smallskip

\item{1)}  A distribution $u\in \cd'(X)$ is {\sl  $\pp$-pluriharmonic}
if $\D_A u=0$ for all $A\in \plp$, or equivalently (see Appendix B)
the $S(\plp)$-Hessian of $u$ is identically zero.
\medskip

The appropriate extensions of partial and strict are more problematic. Uniform strictness can
be put in a satisfactory state.

Suppose $u\in \PSH(X)$ and $\e>0$.  Then $u$ is {\sl $\e$-strict} is either of the following two equivalent conditions are satisfied. ( The proof of this equivalence is omitted.)
\smallskip

\item{a)}  $u-\e|x|^2\in\PSH(X)$.

\smallskip

\item{b)}  For each point $x\in X$ and each function $\vf\in C^2$ near $x$
 which is ``superior'' to $u$ in the sense that $u-\vf$ has a local maximum at $x$, one has
 $\Hess_x \vf -\e I \in \pp$.\smallskip

It is convenient to extend strictness  from $C^2$ functions to general \psh functions as follows.
\smallskip

\item{2)} $u\in\PSH(X)$ is said to be {\sl strict} if $u$ is $\e$-strict for some $\e>0$.
\medskip

The major defect of this definition is best understood by the following example.

\Ex{\II.1} Note that the negation of strictness is no longer the appropriate notion of partially
pluriharmonic.  For the standard Laplacian $\D$ on $\rn$, $u$ is strictly subharmonic if the 
absolutely continuous part of the measure $\D u$ is bounded below a.e. by some $\e>0$. 
Hence, $u$ being subharmonic but not strict does not imply  that $u$ is harmonic.

 In the next section we examine the more difficult notion of $\pp$-partially pluriharmonic functions.

\vskip .3in


\centerline{\bf \UU.  The Dirichlet Problem -- Uniqueness. }
\medskip

Here we consider the \dir problem for functions which are ``$\pp$-partially 
pluriharmonic''.  A full discussion of this concept is given below.
However,  for $C^2$-functions $u$ on $X$ this simply means  that 
$\Hess_x u\in \partial \pp$ for each $x\in X$,
and if, furthermore, $\pp=\pp(G)$ is geometrically defined,  it means that
($u$ is $G$-psh and)  at each $x$, 
$\tr_\x \Hess_x u = 0$  for some $\x\in G$. The main result of this section is the following.

\Theorem{\UU.1.  (Uniqueness for the Dirichlet Problem)} {\sl  Suppose
$\pp$ is an elliptic cone and that $K$ is a compact subset of $\rn$.
If $u_1, u_2 \in C(K)$ are $\pp$-partially  pluriharmonic   on $\Int K$, then }
$$
u_1\ =\ u_2\ \  {\rm on}\  \partial K \qquad \Rightarrow \qquad
u_1\ =\ u_2\ \  {\rm on}\  K 
$$

 In order to formulate our definition for non-$C^2$ functions,
it is useful to study functions $v$ with $-\Hess_x v\notin \Int F$.  These are in some sense (to be made precise) dual to the $\pp$-plurisubharmonic functions.

\Def{\UU.2}  Given a closed subset $F\ss\Symn$, the {\sl \dir dual} is defined to be
$$
\ft \ =\ -(\Int F)\ =\ \sim(-\Int F).
$$
Note that  $$\partial F \ =\  F \cap (-\ft).$$

\Lemma{\UU.3}
$$
B\in \wt{\pp} \quad\iff\quad A+B\in\wt \cp\qquad {\rm for\ all\ \ }A\in \pp
$$
\pf
Since $\Int \pp = \pp+\Int \cp$, we have that
$$
B'\notin \Int \pp \ \iff \ B'-A\notin \Int \cp \fa A\in\pp.
$$
Set $B=-B'$.  Then
$$
B \notin -\Int \pp \ \iff \ B+A\notin -\Int \cp \fa A\in\pp.
$$
\qed

In Appendix A we introduce the class of {\sl subaffine functions $\SA(X)$ on $X$},
and we refer the reader there for a full discussion. 
We mention, however,  that a function $w\in C^2(X)$ is   {\sl subaffine} if for all $x\in X$, 
$\Hess_x w\in   \wt \cp$, i.e., $-\Hess_x w\notin \Int \cp$, i.e., $\Hess_x w$ has at least 
one eigenvalue $\geq0$. The following concept is basic to  uniqueness.

\Def{\UU.4}  A function $v\in\USC(X)$ is said to be of  {\sl type $\wt\pp$ on $X$} if 
$$
A+v\in \SA(X) \fa {\rm quadratic \   functions \ } A\in\pp.
$$
Let $\wt\PSH(X)$ denote the space of all such functions.
\medskip

 If $v\in C^2(X)$, then
$$
v  {\ \rm is\ of\ type\ } \wt\pp
\quad\iff\quad   \Hess_x v \in \wt\pp \fa x\in X.
\eqno{(\UU.1)}
$$
This follows since, as remarked above, $A+v\in \SA(X)$ if and only if $A+\Hess_x v \in \wt\cp$ for all
$x\in X$, which by Lemma \UU.3 is true for all $A\in \cp^+$ if and only if $\Hess_x v\in \wt\pp$.

\Remark{\UU.5}  If $\pp=\pp(G)$ is geometrically defined, then
$$
\wt\pp(G)\ =\ \{B\in\Symn : \tr_\x B \geq0 \ {\rm for\ some\ }\x\in G\}.
$$
To see this first note that
$$
\Int \pp(G)\ =\ \{A\in \Symn : \tr_\x A>0\fa \x\in G\},
$$
that is,
$$
\sim \Int \pp(G)\ =\ \{A\in \Symn : \tr_\x A\leq0 \ {\rm for\ some\ } \x\in G\}.
$$
Now set $B=-A$ and apply the definition of $\wt\pp(G)$.
\medskip

To establish the basic properties of this class it is useful to have alternative
definitions  of type $\wt\pp$ functions.

\Lemma{\UU.6}  {\sl  A function $v$ is of type $\wt\pp$ on $X$
if and only if 
\smallskip
\centerline
{
$u+v\in\SA(X) $ for all $u\in C^2(X)$ which are $\pp$-\psh.
}
\smallskip
Moreover, $v\notin \wt\PSH(X)$ if and only if 
 $\exists \, A\in \pp$, $a$ affine, $x_0\in X$, and $\e>0$ such that
 $$
 \eqalign{
 a+A+v \ &\leq\ -\e|x-x_0|^2\ \ \ \  {\rm for\ }x\ {\rm near\ } x_0  \cr
               &=\ 0\qquad\qquad\qquad {\rm at\ }x=x_0.
               }
  \eqno{(\UU.2)}
 $$
}
\pf  
If $u+v\notin \SA(X)$ with $u\in C^2(X)$ of type $\pp$, then by Lemma \VV.2 there exist
$x_0\in X$, $\e>0$, and $a'$ affine with 
 $$
 \eqalign{
 a'+u+v \ &\leq\ -2\e|x-x_0|^2\ \ \ \  {\rm for\ }x\ {\rm near\ } x_0  \cr
               &=\ \ \ \ 0\qquad\qquad\qquad {\rm at\ }x=x_0.
               }
  \eqno{(\UU.3)}
 $$
Now since $u\in C^2(X)$, we have $A=\half \Hess_{x_0} u \in \pp$. 
Using the Taylor series for $u$ about $x_0$ it is easy to see that 
(\UU.3) implies (\UU.2).   Now (\UU.2) implies that there exists $A\in \pp$ with
$A+v\notin\SA(X)$ (i.e., (\UU.2) implies $v\notin \wt\PSH(X)$).
 The last implication needed 
is trivial from  Definition \UU.4.  Namely, if $v\notin \wt\PSH(X)$,
then$\exists\, u\in C^2(X)$ of type $\pp$ with $u+v\notin \SA(X)$.
 \qed

\Def{\UU.7}  A function $u$ such that $u\in\PSH(X)$ and $-u\in\wt\PSH(X)$
will be called {\sl $\pp$-partially pluriharmonic on $X$}.
\medskip

Note that for such functions $u$, since both $u$ and $-u$ are upper semi-continuous on $X$, 
one has $u\in C(X)$. Furthermore, since $\partial \pp = \pp\cap (-\wt\pp)$, if $u\in C^2(X)$,
then $u$ is $\pp$-partially pluriharmonic if and only if $\Hess_x u\in \partial \pp$ for each 
$x\in X$.

Because of the Maximum Principle in Appendix \VV,  Theorem  \UU.1  is an immediate consequence
of the next result.

\Theorem{\UU.8. (The Subaffine Theorem)}  {\sl  Suppose $\pp$ is an elliptic cone.
If $u\in \PSH(X)$ and $v\in \wt\PSH(X)$, then $u+v\in \SA(X)$.}

\pf  Fact (7) above says that  $u$ is the decreasing limit of smooth functions
$u_j$ which are $\pp$-plurisubharmonic.  By the first part of  Lemma \UU.6, $u_j+v$ is subaffine.
Finally, the decreasing limit of subaffine functions is again subaffine.\qed

\vfill\eject


\centerline{\bf \DD.  The Dirichlet Problem -- Existence. }
\medskip

We now investigate the  existence  
of solutions to the natural Dirichlet problem associated with  $\pp$-\psh functions
on a smoothly bounded domain $\O$. For the existence question, we assume the boundary
$\bo$ is strictly $\pp$-convex, a concept   introduced and 
discussed in detail in \S \JJ.  A principle result, Theorem \JJ.4, states that if $\bo$ is strictly
$\pp$-convex, then there exists 
a smooth, strictly $\pp$-\psh function on a neighborhood of $\overline\O$, which is a defining 
function for $\bo$. It is this result that will be used below, and the reader can, for the moment,
take its conclusion as the working assumption.

As before we assume $\pp$ is an elliptic cone.

\Theorem{\DD.1. (The Dirichlet Problem -- Existence)}  {\sl
Suppose  $\O$ is a bounded domain in $\rn$ with a  strictly  $\pp$-convex boundary.
Given $\vf \in C(\bo)$, the function $u$ on $\overline\O$ defined by taking the upper envelope:
$$
u(x)\ =\ \sup \{v(x) : v\in \PSH(\vf)\}
$$
over the family
$$
\PSH(\vf)\ \equiv  \bigl\{v \ : \ v { \rm \  is\ u.s.c.\ on\ \ }\Ob,\ \  v\bigr|_{\O} \in  \PSH^{\rm u.s.c.}(\O)
\ \ {\rm and\ \ }v\bigr|_{\bo}\leq \vf    \bigr\}
\eqno{(\DD.1)}$$
satisfies:
\medskip

1)\ \ \ $u\in C(\Ob)$,

\medskip

2)\ \ \ $u$ is  $\pp$ partially pluriharmonic on $\O$, 

\medskip

3)\ \ \ $u\bigr|_{\bo} = \vf$  on  $\bo$.}

\medskip
\pf
By the Maximum Principle the family  $\PSH(\vf)$ is uniformly bounded above on $\Ob$ by
$\sup_{\bo} \vf <\infty$.  Hence by 5$'$), the u.s.c. regularization $u^*$ of  the upper envelope $u$
of  $\PSH(\vf)$, belongs to $\PSH^{\rm u.s.c.}(\O)$.
That is
$$
u^*\bigr|_{\O}\in\PSH^{\rm u.s.c.}(\O).
\eqno{(\DD.2)}
$$
Let $h$ denote the unique $\D_A$-harmonic solution to the Dirichlet problem for some mollifying Laplacian $\D_A$.  Then $h\in C(\Ob)$, $u\leq h$ on $\Ob$ and $h=\vf$ on $\bo$.
Hence, $u^*\leq h$ on $\Ob$ so that 
$$
u^*\bigr|_{\bo}  \ \leq \ \vf.
\eqno{(\DD.3)}
$$
This proves
\Prop{\DD.2}  {\sl $u^*\in \PSH(\vf)$ and therefore}
$$
u^*\ =\ u \qquad {\sl on\ \ }\Ob.
\eqno{(\DD.4)}
$$

The following barrier argument is taken from Bremermann [B].

\Lemma{\DD.3}  {\sl  The function $u$ on $\Ob$ is continuous at each point of $\bo$, and 
$u\bigr|_{\bo} =\vf$ on $\bo$.}

 \pf   It suffices to show that 
$$
\liminf_{x\to x_0} u(x)\ \geq \ u(x_0) \qquad {\rm for\ all \ \ }  x_0\in   \bo.
$$
because of (\DD.3) above.
 Fix $x_0\in\bo$ and choose a smooth function $\psi\geq0$ with $\psi(x_0)=0$ and 
 $\psi(x)>0$ for  $x\neq x_0$.  Replacing $\psi$ by a sufficiently small scalar multiple of
 $\psi$ we may assume that $\rho-\psi$ is strictly \psh on $\overline\O$, where $\rho$ is the 
 defining function for $\bo$ given by Theorem \JJ.4. 
  Now for each $\e>0$, there exists $C>>0$ so that the function
 $$
 v(x)\ \equiv\ C(\rho(x)-\psi(x)) +\vf(x_0)-\e
 $$
 satisfies
 $$
 v\ =\ -C\psi+\vf(x_0)-\e\ \leq \ \vf \qquad {\rm on\ \ }    \bo.
 $$
 Thus $v\in \PSH(\vf)$.  Consequently,
 $
 v\leq u
 $
 on $\overline\O$, and so
 $$
 \liminf_{x\to x_0} u(x) \ \geq \  \lim_{x\to x_0} v(x)\ =\    \vf(x_0)-\e. 
 $$
 \qed
 
 We now apply an argument of Walsh [W] to prove interior continuity.
 
 \Prop{\DD.5}   \qquad  $u\in C(\overline\O)$.
 
 \pf 
 Let $\N\d \equiv\{x\in \Ob : \dist(x,\bo)<\d\}$ denote an interior $\d$-neighborhood of the boundary
 $\bo$.
 Suppose $\e>0$ is given.
 By the continuity of  $u$  at points of $\bo$   and the compactness of $\bo$, it follows easily
that there exists a $\d>0$ such  that 
$$
{\rm If \ \ }
x\in \N{2\d}, \ |y|<\d\ {\rm and \ \ } x+y\in \Ob,
\ \ {\rm then\ \ } u(x+y)-u(x)\ <\ \e.
 \eqno{(\DD.5)}
 $$
Now for $|y|<\d$ fixed, consider the function
$$
f_y(x)\equiv \max\{u(x+y)-\e, u(x)\} \qquad{\rm on} \ \ \O-\overline{\N{\d}}.
$$
Note that $f_y\in \PSH(\O-\overline{\N{\d}})$ by 1) in Section \II.

Now consider the restriction of $f_y$ to $\N{2\d}-\overline{\N{\d}}$.  
Then $x\in \N{2\d}$, $|y|<\d$, and $x+y\in\Ob$, so that (\DD.5) implies that
$$
f_y(x)\ =\ u(x)\qquad{\rm on\ \ } \N{2\d}-\overline{\N{\d}}.
$$
We extend $f_y$ to all of $\Ob$ by setting $f_y=u$ on $\N{2\d}$.  The function $f_y$ now belongs to the family $\PSH(\vf)$.  Hence, $f_y\leq u$.  For $x\in\O-\overline{\N{\d}}$ this yields
$$
u(x+y)-\e\ \leq\ u(x)\qquad{\rm if \ } |y|<\d.
$$
Replacing $y$ by $-y$ and x by $x+y$ yields
$$
u(x)-\e\ \leq\ u(x+y)\qquad{\rm if \ } |y|<\d {\rm \  and\ } x\in\O-\overline{\N{2\d}}.
$$
This proves that
$$
|u(x+y)-u(x)|\ <\ \e  \qquad{\rm if \ } |y|<\d {\rm \  and\ } x\in\O-\overline{\N{2\d}}.
$$
\qed\medskip

Finally, to complete the proof of Theorem \DD.1 we must show that the Perron function
$u$ is  $\pp$-partially pluriharmonic on $\O$.
We already have $u\in \PSH(\O)$.  Hence, we must show that $-u \in \wt\PSH(\O)$.
Suppose $K\ss \O$ is compact and let $w$ be a polynomial of degree two which is $\pp$ \psh  
with $w\leq u$ on $\partial K$.  We must show that  $w\leq u$ on $K$. However, this must hold, since otherwise one could change $u$ to $\max\{w,u\}$ on $K$ and violate the  maximality of the Perron function
$u$.\qed

 \Remark{\DD.6}  Suppose $\pp_0\ss\pp_1$ are elliptic cones. Then if a boundary $\bo$ is   strictly
 $\pp_0$-convex, it is also strictly $\pp_1$-convex. 
Furthermore,  if  $u$ is  $\pp_0$-\psh, then it is also $\pp_1$-plurisubharmonic.  It follows that
if $\O\ss\rn$ is a bounded domain with strictly $\pp_0$-convex boundary, and $\vf\in C(\bo)$ is
given, then the unique solutions to the Dirichlet Problem $u_0$ and $u_1$ given by Theorem
\DD.1  for $\pp_0$ and $\pp_1$ respectively satisfy
$$
u_0\ \leq \ u_1\quad{\rm \ \ on\ }\ \O
$$

\vfill\eject

\centerline{\bf  \KK.   $\pp$-Convex Domains}
\medskip

In this section we introduce  the notion of $\pp$-convex domains and give
several characterizations of them.   We then establish 
topological restrictions on any such domain.  In many cases these
restrictions are known to be sharp.  

We assume throughout this section $X$ is a connected open subset of $\rn$, and
that $\pp\subset \Symn$ is a convex cone which satisfies the Positivity Condition:
$\pp\ss\cp$, but  not necessarily the full Ellipticity Condition (i.e., not the Completeness Condition).

\Def{\KK.1}  Given a compact subset $K\ss X$, we define the {\sl $\PSH^{\infty}(X)$-hull of 
$K$} to be the set 
$$
\wh K \ \equiv\ {\wh K}_{\pp, X} \ \equiv\ \{x\in X : u(x) \leq \sup_K u \ \ {\rm for\ all\  } u\in\PSH^\infty(X)\}.
$$
If $\wh K=K$, then $K$ is called {\sl $\pp$-convex}. 

\Lemma {\KK.2}  {\sl Suppose $K$ is a compact subset of $X$.  A point $x$ is not in $\wh K$
if and only if there exists $u\in \PSH^\infty(X)$ with $u\geq 0$ on $X$ and 
$u=0$ on a neighborhood of $K$ but $u(x)>>0$; and with $u$ strict at $x$.}

\pf Suppose $x_0\notin \wh K$.  Then there exits $v\in \PSH^\infty(X)$ with $\sup_K v <0<v(x_0)$.
Multiplying $v$ by a large constant, we may assume that $v(x_0)$ is large.
Replacing $v$ by $v+\e|x|^2$, we may assume that $v$ is strict at $x_0$.
An $\e$-approximation $u=\max_\e\{0,v\}$ to the maximum $\max\{0,v\}$ satisfies
all the conditions.\qed

\Prop{\KK.3}  {\sl  The following two conditions are equivalent.\smallskip

\item {1)} $K\ss \ss X \ \ \Rightarrow\ \ \wh K\ss\ss X$.  \smallskip

\item {2)} There exists a $C^\infty$ proper exhaustion function $u$ for $X$ which is
strictly $\pp$-psh.  \smallskip
}

\Def{\KK.4}  If the equivalent conditions of Proposition \KK.3 are satisfied, then $X$ is a {\sl
$\pp$-convex domain } in $\rn$.

\medskip
\noindent
{\bf Proof of Proposition \KK.3}.  
We first show that $2)\Rightarrow 1)$.  If $K\ss X$ is compact, then $c=\sup_K u$ is finite and $\wh K$ is contained in the compact prelevel set $\{u\leq c\}$.

To see that  $1)\Rightarrow 2)$, choose compact $\PSH^\infty(X)$-convex sets
$K_1\ss K_2\ss \cdots$ with $K_m\ss K_{m+1}^o$  and   $X=\bigcup_m K_m$.
By the Lemma above and the compactness of $K_{m+2}- K_{m+1}^o$ we may find
$u^1,...,u^N\in \PSH^\infty(X)$, which are non-negative and vanish on a neighborhood of $K_m$,
with $u_m = \max_\e\{u^1,...,u^N\} > m$ on $K_{m+2}- K_{m+1}^o$.  The maximum
$u=\max\{u_1,u_2,...\}$   satisfies 2), except for strictness. To obtain strictness, replace
$u$ by $u+\half |x|^2$, which is strict because $I$ is an interior point of $\cp\ss\pp$.  \qed

\Remark {\KK.5}  Condition 2) in Proposition \KK.3 can be weakened in several ways.

First, the exhaustion $u$ need only be $\pp$-plurisubharmonic, not strict,  
since one can always replace $u$ with   $u+  |x|^2$.

Second, 
$u$ only needs to be defined near $\infty$ in the one point compactification of $X$.
More precisely, if there exists $u\in \PSH^{\infty}(X-K)$, where $K$ is compact, $u$
is bounded near  $K$,  and  $\lim_{x\to\infty} u(x) = \infty$, then 2) holds.
To see this, note that for large $c$, $v=u+|x|^2$ 
is a smooth strictly $\pp$\psh function outside the compact subset
$\{v\leq c-1\}$.  Pick a convex increasing function $\vf\in C^\infty(\bbr)$ with $\vf = c$ on a neighborhood of $(-\infty, c-1]$ and with $\vf(t)=t$ on $(c+1,\infty)$.  Then $\vf(v(x))\in \PSH^\infty(X)$
and equals $v(x)$ outside the compact set $\{v\leq c+1\}$.\medskip

\vskip.3in

\centerline{\bf  \RR.   Topological Restrictions on $\pp$-Convex Domains}
\medskip

We begin our discussion of the topology of $\pp$-convex domains with the following definition.
Note that for any linear subspace, $W\ss \rn$ there is a natural embedding $\Sym(W)\ss\Symn$.

\Def{\RR.1}  

a) A linear subspace $W\ss\rn$ is {\sl $\plp$-free} if $\plp\cap \Sym(W) =\{0\}$.
In the geometric case  where $\plp=\plp(G)$, this means that $W$ does not contain any $p$-planes 
$\xi\in G$.  In this case we say that $W$ is {\sl $G$-free}.\smallskip

b) A linear subspace $N\ss\rn$ is {\sl $\pp$-strict} if $P_N\in\Int\pp$.

\Lemma{\RR.2}  {\sl  Suppose that $\rn = N\oplus W$ is an orthogonal decomposition.
Then $W$ is $\plp$-free if and only if $N$ is $\pp$-strict.
\smallskip
}

\pf
If $N$ is not strict, then by the Positivity Condition
 $P_N \in\partial \pp$.
Hence, there exists $A\in\plp$, $A\neq 0$, with $\langle P_N, A \rangle=0$. 
 By the positivity assumption $\plp\ss\cp$  and the basic fact  (\EE.3),
 it follows easily that  $\langle P_N, A \rangle=0$
if and only if $A\in \Sym(W)$.  Thus, $\plp\cap \Sym(W)\neq \{0\}$, contradicting $W$ being free.
On the other hand, if $P_N$ is strict, then for all $A\in\plp$, $\langle P_N, A\rangle >0$ unless $A=0$, proving that
$\plp\cap \Sym(W) =\{0\}$. \qed

\Remark{\RR.3}  
\smallskip
\centerline
{
$P_N$ is strict if and only if $\Int\pp\cap\Sym(N) \neq \emptyset$.
}
\pf
  Note that if 
 $P_N$ is strict, then $P_N\in\Int\pp\cap\Sym(N)$. For the converse, suppose
there exists $H\in \Int\pp\cap\Sym(N)$, then $H\neq 0$ and $\langle H, A\rangle >0$
for all non-zero $A\in\plp$.  However,  $\langle H, A\rangle =0$ for all $A\in \Sym(W)$
proving that $W$ is free. Hence  $N$ is strict by Lemma \RR.2.

\Def{\RR.4}  The {\sl free dimension of  $\plp$}, denoted by  \fd$(\plp)$
(or \fd$(G)$ in the geometric case),  is the maximal dimension of a $\plp$-free
subspace of $\rn$. By Lemma \RR.2 this equals the maximal codimension
of a $\pp$-strict subspace.\medskip

Somewhat surprisingly the Andreotti-Frankel Theorem in complex analysis 
has a very general extension.
 The usual proof of the Andreotti-Frankel result is quite specific to complex analysis, relying
 on   canonical forms.

\Theorem{\RR.5}  {\sl Let $X$ be a $\pp$-convex domain in $\bbr^n$.  
Then $X$ has the homotopy-type of a CW-complex of dimension $\leq $  \fd$(\plp)$.
}

\pf
Let $u\in \cix$ be a   strictly $\pp$-\psh proper exhaustion function. By standard approximation 
theorems (cf. [MS]) we may assume that all critical points of $u$ are non-degenerate.
The theorem will follow if we can show that each critical point $x_0$ of $u$ in $X$  has index 
$\leq $ \fd$(\pp)$.

Since $u$ is $\pp$-\psh, we have $\Hess_{x_0} u\in \pp$, that is 
$$
\langle \Hess_{x_0} u, A\rangle\ \geq\ 0  \fa A\in \plp.
\eqno{(\RR.1)}
$$
Suppose now that  the index of $\Hess_{x_0} u$ is $>$ \fd$(\plp)$. Then there exists a subspace
$W\ss \rn$  with $\dim(W)> $ \fd$(\plp)$  such that
$$
\Hess_{x_0} \left( u\bigr|_W\right)\  <\ 0.
\eqno{(\RR.2)}
$$ 
Now by definition
of \fd$(\plp)$ there exists a non-zero $A\in \Sym(W)\cap \plp$.  
Hence, 
$\langle A,  \Hess_{x_0} u\rangle  = 
\left\langle A,  \Hess_{x_0} \left( u\bigr|_W\right)\right\rangle  
< 0 $. \qed

\Remark{\RR.6}  The free dimension of $\plp$ is $n-1$ if and only if each hyperplane
$W$ is free, i.e., $\plp\cap \Sym(W)=\{0\}$, or equivalently each $P_e\in\Int\pp$ for 
$0\neq e\in\rn$. Otherwise the free dimension is $< n-1$.  In this case  $\bo $ is connected 
for every $\pp$-convex domain. (This is the case $k=0$ in the next Corollary.)  
A special case of this connectedness appears as Lemma A in [CNS].

\Cor{\RR.7} {\sl Let $\O\ss\ss X$ be a strictly \ppc\ domain with smooth boundary $\bo$, and 
let $D$ be the free dimension of $\plp$.  Then 
$$
H_k(\bo  ; \,\bbz)\ \cong \ H_k(\O  ; \,\bbz)\qquad  {\rm for\ all\ } k< n-D-1
$$
and the map $H_{n-D-1}(\bo;\, \bbz)\to H_{n-D-1}(\O;\, \bbz)$ is surjective.}

\pf This follows from the exact sequence
$$
H_{k+1}(\O, \bo  ; \,\bbz)\ \to \ H_k(\bo  ; \,\bbz)\ \to \  H_k(\O  ; \,\bbz)\ \to \ H_k(\O,\bo  ; \,\bbz),
$$
Lefschetz Duality: $H_k(\O,\bo;\,\bbz) \cong H^{n-k}(\O;\,\bbz)$, and Theorem \RR.5.\qed

\bigskip
\noindent
{\bf Geometric Examples.} Consider the geometric case $\pp=\pp(G)$.  Set    fd$(G)$ = \fd$(G)$.
The following facts were shown in [HL$_{2}$].
\medskip

\noindent
\item{1.}\  $G=G(1,\rn)$ (Convex geometry).\     fd$(G)  =0$.

\medskip

\noindent
\item{2.}\  $G=G(n,\rn)$ ($\PSH(X,G)$ = {subharmonic  functions} on $X$). \   fd$(G) =  n-1$.

\medskip

\noindent
\item{3.}\  $G=G(p,\rn)$ for $1<p<n$.   \  fd$(G)  = p-1$.
\medskip

\noindent
\item{4.}\  $G= \bbp^{n-1}(\bbc) = G_\bbc(1,\bbc^n) \subset G(2, \bbr^{2n})$ (Complex psh-functions).
  \  fd$(G)  = n$.\medskip

\noindent
\item{5.}\  $G= \bbp^{n-1}(\bbh) = G_\bbh(1,\bbh^n) \subset G(4, \bbr^{4n})$ (Quaternionic   psh-functions).   \  fd$(G)  = 3n$.
\medskip

\noindent
\item{6.}\  $G = G_\bbc(p,\bbc^n)$ for $1<p<n$.
 fd$(G) =  n+p-1$.
\medskip

\noindent
\item{7.}\  $G = G_\bbh(p,\bbh^n)$ for $1<p<n$.
 fd$(G) =  3n+p-1$.
\medskip

\noindent
\item{8.}\  $G = \{x_1$-axis$\} \subset G(1,\bbr^n)$.  fd$(G) =  n-1$.
\medskip

\noindent
\item{9.}\  $G = {\rm SLAG} \subset G(n,\bbc^n)$,  the  special Lagrangian 
 $n$-planes in $\bbc^n$.  fd$(G) =  2n-2$
 \medskip

\noindent
\item{10.}\  $G = {\rm ASSOC} \subset G(3,\bbr^7)$,  the  associative 3-planes in 
Im$\bbo \cong \bbr^7$.  fd$(G) =  4$.
 \medskip

\noindent
\item{11.}\ $G = {\rm COASSOC} \subset G(4,\bbr^7)$,  the coassociative 4-planes in 
Im$\bbo$.  fd$(G) =  4$.
 \medskip

\noindent
\item{12.}\  $G = {\rm CAY} \subset G(4,\bbr^8)$,  the  Cayley 4-planes in 
the octonions  $\bbo\cong\bbr^8$.  fd$(G) =  4$.
 \medskip

\noindent
\item{13.}\  $G = {\rm LAG} \subset G(n,\bbc^n)$,  the set of   Lagrangian 
 $n$-planes in $\bbc^n$.  fd$(G) =  2n-2$.
 \medskip

\bigskip
\noindent
{\bf Some Non-Geometric Examples.}  
Let $\s_k(A):\Symn\to \bbr$ be the $k$th elementary symmetric
function of the eigenvalues defined by the equation $\det(I+tA)=\sum_k \s_k(A) t^k$
Consider the closed cone $\pp(\s_k)$ whose interior
is the connected component, containing $I$, of the set
$
 \{A\in\Symn :  \s_k(A)>0\}
$. Then \medskip
\item{14.}\   fd$(\pp(\s_k)) = n-k$.
\medskip

Note that every $k$-plane $N$ is $\pp(\s_k)$-strict because $\s_k(P_N)=1$.
On the other hand $\s_k(P_N)=0$ for any $(k-1)$-plane. 
Thus, the strict dimension of $\pp(\s_k)$  is $k$ or equivalently, the free dimension
of $\plp(\s_k)$ is $n-k$.

%
%

\vskip .3in

\centerline{\bf  \LL.   $\plp$-Free Submanifolds}
\medskip

We have seen in \S \RR \ that there are sometimes quite strong   restrictions
on the homotopy dimension of $\pp$-convex domains. 
In this section we show that within these restrictions the topology of such 
domains can be quite complicated. One of the main results, Theorem \LL.4,  is that any
submanifold $M \ss X$, which is $\plp$-free,  has a fundamental system of strictly
$\pp$-convex neighborhoods  homotopy equivalent to $M$.
 
 Most  proofs in this section are omitted since they carry over by direct generalization from [HL$_2$]. 
The reader can consult  [HL$_2$] for further results and details.

\Def{\LL.1}  A closed submanifold $M\ss X\ss \rn$ is {\sl $\plp$-free} if the tangent space
$T_xM$ is $\plp$-free at each point $x\in M$.
(In the geometric case where $\plp=\plp(G)$ this means that there are no $G$-planes which
are tangential to $M$.)

\Theorem{\LL.2}  {\sl Suppose $M$ is a closed submanifold of $X\ss\rn$.  Then $M$ is 
$\plp$-free if and only if the square of the distance to $M$  is strictly $\pp$-\psh at each point in $M$  (and hence in a neighborhood of $M$ in $X$).}

\pf
Given $x_0\in M$, let $N=(T_{x_0}M)^\perp$ denote the normal to $M$ at $x_0$.  Let $f_M(x) =
\half \dist^2_M(x)$ denote half  the square of the distance to $M$.  One can calculate that
$$
\Hess_{x_0} f_M \ =\ P_N.
$$
(See, [HL$_2$, (6.3)].) Now the theorem is immediate from  Lemma \RR.2.\qed

\Theorem{\LL.3}  {\sl Consider the two classes of closed sets.
\smallskip

\item{1)}  Closed subset $Z\ss M$ of a $\plp$-free submanifold $M\ss X$.

\smallskip

\item{2)}  Zero sets $Z=\{f=0\}$ of non-negative strictly  $\pp$-\psh functions $f$. \smallskip

\noindent
Locally these two classes are the same.}

\pf
Suppose $Z\ss M$ is as in  1).  Choose $\psi\in\cix$  with $\psi\geq 0$  and $\{\psi=0\}=Z$. Then  
for $\e>0$ small, the function $f_M+\e\psi$ is strictly $\pp$-\psh and $Z=\{f_M+\e\psi=0\}$.

Assume $Z=\{f=0\}$ is as in 2).  At $x_0\in Z$ choose coordinates $x=(z,y)$ in a 
neighborhood of $x_0$ so that 
$$
\Hess_{x_0} f \ =\ \left(\matrix{0&0\cr 0&\Lambda}\right)
$$
where $\Lambda$ is a diagonal matrix with non-zero diagonal entries:
$$
{\partial^2 f\over \partial y_1^2}(x_0),\ ...\ ,{\partial^2 f\over \partial y_r^2}(x_0).
$$
Set 
$$
M\ =\ \left\{ {\partial f\over \partial y_1} = \cdots = {\partial f\over \partial y_r} =0\right\}.
$$
This defines a submanifold $M$ in a neighborhood of $x_0$, 
since $\nabla {\partial f\over \partial y_1} , ..., \nabla {\partial f\over \partial y_r} $ are linearly independent
at $x_0$.  At $x_0$ the normal space to $M$ is 
$
N\ =\ \{(0,y) : y\in \bbr^r\}.
$
Strict plurisubharmonicity implies  $\Hess_{x_0} f \in \Int \pp\cap \Sym(N)$ 
so that $T_{x_0}M = N^\perp$ is $\plp$-free by Lemma \RR.2. Since the freeness  condition is open, the manifold $M$ is $\plp$-free in a neighborhood of $x_0$. Since $f\geq 0$, $\nabla f =0$ at all points of $Z=\{f=0\}$, and so $Z\ss M$.
 \qed

\Theorem{\LL.4}  {\sl Suppose $M$ is a $\plp$-free closed submanifold of $X\ss \rn$.  
Then there exists a fundamental neighborhood system $\cf(M)$ of $M$ 
consisting of $\pp$-convex domains.  Moreover,
\smallskip
\item{a)}  $M$ is a deformation retract of each $U\in \cf(M)$.

\smallskip
\item{b)} Each compact subset $K\ss M$ is $\PSH^\infty(U,\pp)$-convex for  each $U\in \cf(M)$.
}

\medskip
The proof of this theorem is exactly as in [HL$_2$, Thm. 6.6]  and is omitted.


\vfill\eject

\centerline{\bf \JJ.    $\pp$-Convex Boundaries   }
\medskip

In this section we introduce the notion  of    $\pp$-convexity for smooth
boundaries of domains in $\rn$.   We show, for  bounded domains,  that if the boundary 
is strictly \ppc\ at each point, then there exists a global defining function
$\rho$ for the domain which is  strictly $\pp$-\psh on its closure.  It is then easy to see
that $-\log(-\rho)$ is a strictly $\pp$-\psh exhaustion, and so the domain is \ppc.

   Fix a domain $\O\subset\subset\rn$ with smooth boundary $\bo$.
   By a   {\sl  defining function for} $\bo$ we mean a smooth  function $\rho$ defined in a neighborhood of $\bo$ such that in this neighborhood
 $\O=\{x\in\rn: \rho(x)<0\}$ and
   $\n \rho \neq 0$ on $\bo$.
  
  An element $A\in \plp$ is said to be {\sl tangential at $x\in\bo$ } if $A\in \Sym(T_x\bo)$.
  In terms of the  $2\times2$  blocking induced by the decomposition 
  $\rn= N_x(\bo)\oplus T_x(\bo)$, this means $A=\left(\matrix {0& 0 \cr 0&a\cr}\right)$

 \Def{\JJ.1} We say that $\bo$ is  strictly {\sl $\pp$-convex}  at a point $x\in\bo$ if
 $$
 \langle \Hess_x \rho , A\rangle \  > \ 0 \qquad{\rm for\ all\ non\,zero\ } A \in\plp \ {\rm which\  are\  tangential\  at\ }x.
\eqno{(\JJ.1)} 
 $$
 If $ \langle \Hess_x \rho , A\rangle  \geq 0$ for all  tangential  $A\in\plp$, then 
 $\bo$ is said to be {\sl   $\pp$-convex at $x$}.\medskip

\Note{\JJ.2} These notions are independent of choice of defining function $\rho$. If $\wt \rho = f\rho$ with 
$f>0$ in $C^\infty(\bo)$, then $\Hess {\wt\rho} = f \Hess \rho +\rho \Hess f + 2 \nabla \rho\circ\nabla f$,
and so $ \langle \Hess_x {\wt\rho} , A\rangle = f \langle \Hess_x \rho , A\rangle$ for $A\in\plp$
which are tangential at $x$.

 \Remark{\JJ.3}  ({\bf The Geometric Case})  If   $\plp=\plp(G)$, where $G$ is a closed subset
 of the Grassmannian $G(p,\rn)$, note that $A\in \plp(G)$ is tangential if and only if 
 $A=\sum_j\l_j P_{\x_j}$ with each $\l_j>0$ and each 
 $\x_j\in G$     {\sl tangential} in the sense that $\span \x_j \ss T_x\bo$.
 To show this let $n$ denote a unit normal to  $\bo$ at $x$.  If $A\in \plp(G)$, then by definition
 $A=\sum_j\l_j P_{\x_j}$ with each $\l_j>0$ and each 
 $\x_j\in G$.  If $A$ is tangential to $\bo$ at $x$, then $0=\langle A, P_n\rangle =
 \sum_j\l_j \langle P_{\x_j}, P_n\rangle$ and hence each $\langle P_{\x_j}, P_n\rangle = |n\hk \x_j|^2$
 vanishes, which implies that $\span \x_j \ss T_x\bo$.  Consequently,
{\sl $\bo$ is  strictly \ppc\ at $x\in\bo$ if and only if 
$$
\tr_\x\Hess_x \rho \ =\ \langle  \Hess_x\rho, P_\x\rangle\  > \ 0 \fa \x\in G \ {\rm which\ are\ tangential\ at\ } x
\eqno{(\JJ.2)} 
 $$
 (and $\bo$ is  \ppc\ at $x$ iff  $\tr_\x\Hess_x \rho \geq 0$ for all  $\x\in G$ tangential at $x$).}

 \Theorem {\JJ.4} {\sl Suppose  that $\O$ has a   strictly $\pp$-convex boundary. 
 Then there exists a strictly $\pp$-\psh function on a neighborhood of $\overline\O$ which is a defining function for $\bo$.}
 
 \pf   
Fix $C>0$ and consider $\wt \rho = \rho+\half C\rho^2$. This is also a defining function for 
$\bo$.  At $x\in \bo$
$$
\Hess_x \wt \rho \ =\   \Hess_x \rho +C(\n\rho\circ\n\rho).
\eqno{(\JJ.3)}
$$

\Lemma{\JJ.5}  {\sl For $C$ sufficiently large, $\wt \rho$ is strictly $\pp$-\psh at each point $x\in \bo$.}

\pf
Since $\plp\ss\cp$, condition (\EE.3) states that the tangential condition 
$$
A\in \plp\cap \Sym(T_x\bo)\quad{\rm is\  equivalent\ to\ }\quad
A\in\plp\ \ {\rm and \ \ } \langle \n\rho(x)\circ\n\rho(x), A \rangle =0.
\eqno{(\JJ.4)}
$$
Now restrict attention to the compact base 
$\cb_+ = \{A\in \plp : \tr A=1\}$ for $\plp$.
Consider the open subsets of  $\bo\times \cb_+$ defined by
$$
U_\delta \ =\ \{(x,A)\in \bo\times \cb_+ : \langle \n\rho(x)\circ\n\rho(x), A \rangle <\delta\}.
\eqno{(\JJ.5)}
$$
   Because of (\JJ.4) these sets $U_\delta$ form a fundamental
 neighborhood system, in $\bo\times \cb_+$,  for the compact set
 $$
 K\ =\ \{(x,A) \in  \bo\times \cb_+  : A\ {\rm is\ tangential\ to \ }\bo \ {\rm at\ }x\}.
 $$
The hypothesis that $\bo$ is strictly $\plp$-convex implies that, for
$\e>0$ sufficiently small, $N(K) =\{(x,A) \in  \bo\times \cb_+  :  \langle \Hess_x\rho, A\rangle >\e\}$
contains $K$.
This proves that there exist $\e, \delta>0$ such that for each 
$(x,A) \in  \bo\times \cb_+ $
$$
\langle \n\rho(x)\circ\n\rho(x), A \rangle <\delta\quad\Rightarrow\quad 
 \langle \Hess_x\rho, A\rangle >\e.
\eqno{(\JJ.6)}
$$
Choose $M>0$ so that $-M < \langle \Hess_x\rho, A\rangle$ for all $(x,A) \in  \bo\times \cb_+ $.
Then, for $\langle \n\rho(x)\circ\n\rho(x), A \rangle \geq \delta$, one has
$$
   \langle \Hess_x  \wt\rho, A\rangle \ =\ 
  \langle \Hess_x\rho +C( \n\rho(x)\circ\n\rho(x))      , A\rangle   \ \geq\ C\delta-M,
\eqno{(\JJ.7)}
$$
while for $\langle \n\rho(x)\circ\n\rho(x), A \rangle < \delta$, one has
$$
   \langle \Hess_x  \wt\rho, A\rangle \ \geq\ 
  \langle \Hess_x\rho      , A\rangle   \ >\  \e
\eqno{(\JJ.8)}
$$
by (\JJ.6).
Choose $C>M/\delta$.
\qed\medskip

 Since $\wt\rho$ is strictly \ppsh  at each point $x\in\bo$, the same is true in a neighborhood
 $ \{-2t<\wt\rho<2t\}$        of $\bo$.

 To complete the proof of the theorem, 
 it remains to extend 
$\wt\rho$ to all of $\overline\O$.  The function $\max\{\wt\rho, -t\}$ is a 
$\pp$-plurisubharmonic extension, but it is not smooth or strict.
However, replacing $-t$ by $a|x|^2-t$, where $a>0$ is chosen small enough so that 
$a|x|^2-t < \wt\rho$ on $\{-{t\over 2}<\wt\rho<0\}$, and then smoothing, we have that for 
$\e>0$ sufficiently small,
$$
\wh \rho \ =\ \max_\e\{\wt\rho, a|x|^2-t\}
$$
is a $C^\infty$ strictly $\pp$-\psh function on a neighborhood of $\overline \O$ which agrees with 
$\wt\rho$ on a neighborhood of $\bo$.\qed

\Remark{\JJ.6} In the non-geometric cases, where $\pp$ is given but $\plp$ may be difficult to determine explicitly, the proof of Lemma \JJ.5 (see (\JJ.3)) provides a convenient criterion for strict boundary convexity.  Namely: 
$$
\eqalign{
&\bo \ \ {\rm is\ strictly\ } \pp \ {\rm convex\ at\ } x\in \bo \qquad\Leftrightarrow  \cr
 \Hess_x \rho\, +\,&C(\n \rho(x) \circ \n \rho(x)) \in \Int \pp\quad {\rm for\ } C>0 \ {\rm sufficiently \ large} \cr
}
\eqno{(\JJ.9)}
$$
The corresponding statement for $\pp$-convexity is false.  Consider $\pp=\cp$ and $n=2$ with 
$T_x\bo = \span e_2$.  Then $H=\left( \matrix{ 0&a\cr a&0\cr}\right)$ is $\geq 0$ and tangential
at $x$, but $H+ C e_1\circ e_1 = \left( \matrix {C&a\cr a&0\cr}\right)$ is never in $\pp=\cp$.

\medskip

We now consider convexity of $\bo$  in terms of its  second fundamental form $II$
with respect to the outward pointing normal.  Let $\rho$ denote the signed distance function
to $\bo$, i.e., $\rho(x) =  -\dist(x, \bo) $ for $x\in \O$ and  $\rho(x) =  \dist(x, \bo) $ for $x\notin \O$
(so $\rho$ is a defining function for $\bo$).
One computes (see [HL$_2$, \S 5])  that for points   $x\in\bo$ 
$$
\Hess_x \,  \rho \ =\   \left( \matrix{ 0&0\cr 0&  -II  \cr}\right)
\eqno{(\JJ.10)}
$$
with respect  to the orthogonal decomposition
$$
T_x\rn \ =\ N_x\bo \oplus T_x\bo.
\eqno{(\JJ.11)}
$$

\Prop{\JJ.7}  {\sl
Suppose     $\O\ss\ss \rn$ has a smooth boundary, and denote by $II$
the second fundamental form of $\bo$ with respect to the outward pointing normal $n= \n \rho$.
Then $\bo$ is strictly $\pp$-convex at a point $x\in \bo$ if and only if 
$$
\langle  II, A \rangle \ < \ 0 \fa {\sl non\, zero\ } A \in \plp \ {\sl which\   are\  tangential\ at\  } x
$$
or equivalently
$$
 -II  +  C n \circ n \in \Int \pp \qquad{\rm for}\  C >0\ \ {\sl sufficiently\ large}.
\eqno{(\JJ.12)}
$$
}

\pf Since $\rho$ is a defining function for $\bo$,  the first assertion follows immediately from  (\JJ.10).
The proof of (\JJ.12) is discussed in Remark \JJ.6.  \qed

\Remark {\JJ.8. (The Geometric Case)}  The boundary $\bo$ is strictly $\pp(G)$-convex at $x\in\bo$
 if and only if 
$$
\tr_\x II\ <\ 0 \quad{\rm for\ each\ } \x\in G \ {\rm which \  is \  tangential \ at\ } x.
$$

\medskip

Finally we discuss the relationship of boundary convexity to the convexity of the domain itself.

\Prop{\JJ.9}  {\sl Suppose that $\O\subset\ss \rn$ has a smooth, strictly \ppc\ boundary.
Then $\O$ is a \ppc\ domain.}

\pf If $\rho$ is a strictly \ppc\  defining function (such as the one given by Theorem \JJ.4),
then $-\log (\delta)$, with $\delta = -\rho$,  is a strictly \ppsh exhaustion function. One computes that
$$
\langle\Hess(-\log \delta), A\rangle  \ =\ {1\over \delta}\langle\Hess\, \rho, A\rangle 
+ {1\over \delta^2}\langle\n\rho \circ \n \rho, A\rangle .
\eqno{(\JJ.13)}
$$
The right hand side is $>0$ for all non-zero $A\in \plp$.\qed
\medskip

In general it is not true that boundaries of $\pp$-convex domains are $\pp$-convex.
(See [HL$_2$, \S 5]  for examples).  However, we have the following.

\Theorem{\JJ.10} {\sl
Let  $\delta = \dist(\bullet, \bo)$ denote the   distance to $\bo$ in $\O$. If $-\log \d$ is $\pp$-\psh
near $\bo$, then $\bo$ is $\pp$-convex.}

\pf
If $\bo$ is not $\pp$-convex, then there exists $x\in\bo$ and $A\in \plp\cap \Sym(T_x\bo)$
 with 
$ \langle  II, A \rangle \ >\ 0$. 
Since $A$ is tangential, we have 
$\langle\nabla \delta \circ \nabla \delta, A  \rangle =0$. 
Let $\ell$ denote the line segment in $\O$ which emanates from $x$ normally to the boundary,
i.e., in the direction $\n\delta$.
It   follows from (\JJ.10) and (\JJ.13)  that 
$\langle \Hess (-\log \delta), A  \rangle 
= - {1\over \delta} \langle \Hess\, \d, A  \rangle
<0$  at all points of $\ell$ near to $x$.
Consequently, $-\log \delta$ is not $\pp$-\psh  in any neighborhood of $\bo$.\qed


\vfill\eject

\centerline{\bf  Appendix \VV. }
\smallskip
\centerline{\bf   The Maximum Principle and Subaffine Functions. }
\medskip

 An upper semi-continuous function 
$u:X\to [-\infty, \infty)$ satisfies the {\sl maximum principle} if for each 
compact subset $K\ss X$
$$
\sup_K u\ \leq \ \sup_{\partial K} u.
\eqno{(\VV.1)}
$$
 A function $u$ may locally satisfy the maximum principle
  without satisfying the maximum principle on all of $X$.
(Consider, for example, a function $u$ on $\bbr$ with compact support, $0\leq u\leq 1$,
$u\equiv 1$ on a neighborhood of the origin and otherwise monotone.)  However, this situation is easily remedied.  First note that (\VV.1) is equivalent 
to the condition that:
$$
u\ \leq \  c \ \ {\rm on\ } \partial K \qquad\Rightarrow\qquad 
u\ \leq \  c \ \ {\rm on\ } K  \qquad
 {\rm for\ all\ constants\ } c,
\eqno{(\VV.1)'}
$$
 i.e., $u$ is {\sl sub-constants}.   Replacing the constant functions by the affine functions, consider the condition:
 $$
u\ \leq \  a \ \ {\rm on\ } \partial K \qquad\Rightarrow\qquad 
u\ \leq \  a \ \ {\rm on\ } K \qquad
 {\rm for\ all\ affine\ functions\ } a
\eqno{(\VV.2)}
$$
\Def{\VV.1}   A function $u\in \USC(X)$ satisfying (\VV.2)
 for all compact subsets $K\ss X$  will be called {\sl subaffine} on $X$. 
Let $\SA(X)$ denote the space of all $u\in \USC(X)$ that are locally subaffine on $X$,
i.e., for all $x\in X$ there exists a neighborhood $B$ of $x$ with $u\bigr|_B$ sub-affine.
\medskip

Note that if $u$ is subaffine on $X$, then the restriction to any open subset is also subaffine.

\Lemma{\VV.2}  {\sl If $u$ is locally subaffine on $X$, then $u$ is subaffine on $X$.
In fact, $u$ is not subaffine on $X$ if and only if }
$$\eqalign{
{\sl There\ exist\ \ }  x_0\in X, \ \ &a\ {\sl affine, \ and\ } \e>0 \ {\sl such\ that\ } \cr
 (u-a)(x)\ \ &\leq\  -\e|x-x_0|^2\ \ {\sl near\ } x_0, \ \ {\sl and}\cr
 (u-a)(x_0)\   &=  \  0\  \cr
 }
\eqno{(\VV.3)}
$$
\pf
Subaffine implies locally subaffine, which implies (\VV.3) is impossible.  
Hence, it remains to show that if (\VV.3) is false, then $u$ is subaffine, or
equivalently,
 if $u$ is not subaffine on $X$, then (\VV.3) is true.
If $u$ is not subaffine on $X$, then for some compact set $K\ss X$ and some affine function
$b$, the difference $w=u-b$ has an interior maximum point for $K$.  For $\e>0$ sufficiently small,
the same is true for $w = u+ \e |x|^2-b$.    Choose a maximum point $x_0\in\Int K$ for $w$
and let $M=w(x_0)$ denote the maximum value on $K$.    
Then $u+ \e |x|^2-b-M \leq 0$ on $K$ and equals zero at $x_0$.  Since $\e|x|^2$ and $\e|x-x_0|^2$
differ by the affine function, this proves that there is an affine function $a$ such that
$u+\e|x-x_0|^2-a\leq0$ on $K$ and is equal to zero at $x_0$, i.e., 
(\VV.3) is true. \qed

\Theorem{\VV.3. (Maximum Principle)}
{\sl  Suppose $K\ss\rn$ is compact and $u\in  \USC(K)$.  If 
$u\in\SA(\Int K)$, then}
$$
\sup_K u\ \leq\ \sup_{\partial K} u.
$$

\pf  Exhaust $\Int K$ by compact sets $K_\e$. Since $u\in \SA(\Int K)$, 
$\sup_{K_\e} u\leq \sup_{\partial K_\e} u$. Since $u\in \USC(K)$, 
each $U_\d = \{x\in K: u(x) < \sup_{\partial K} u + \d\}$,
for $\d>0$, is an open neighborhood of $\partial K$ in $K$.
Therefore, there exits $\e>0$ with $\partial K_\e\ss U_\d$ which implies that 
$\sup_{\partial K_\e} u \leq \sup_{\partial K} u+\d$.  \qed

\medskip

For functions which are $C^2$ (twice continuously differentiable), the subaffine condition is a condition on the hessian of $u$ at each point.

\Prop{\VV.4}  {\sl Suppose $u\in C^2(X)$.  Then }
$$
u\in\SA(X) \ \ \ \iff \ \ \   \Hess_x u \ \ {\sl has \ at \ least \ one \ 
eigenvalue\ }\geq0 \ {\sl at\ each \ point\ }x\in X.
$$
\pf  Suppose $ \Hess_{X_0} u <0$ (negative definite) at  some point $x_0\in X$.
Then the Taylor expansion of $u$ about $x_0$ implies (\VV.3)
Therefore, since $u(x_0)=0$, $u\notin  \SA(X)$.

Conversely, if $u\notin \SA(X)$, then (\VV.3) is valid for some point $x_0\in X$ which implies
that $\Hess_{x_0} u + \e I \leq 0$.     So  $\Hess_x u <0$ is negative definite.
\qed

\medskip

\Ex{(n=1)}  Suppose $I$ is an open interval in $\bbr$.  Then
$$
u\in \SA(I) \quad\iff\quad {\rm either \ \ } u\in {\rm Convex}(I)\ {\rm or\ } u\equiv-\infty.
$$
\pf 
Suppose $u\in\SA(I)$ equals $-\infty$ at one point $\a\in I$ but $u$ is finite at another point
$\b\in I$. Choose $a$ to be the affine function with $a(\a) = -N$ and $a(\b) = u(\b)$.
By (\VV.2),  we have $u\leq a$ on $[\a,\b]$, which implies (by letting $N\to\infty$) that
$u\equiv -\infty$ on $[\a,\b)$.  Hence $u$ is either $\equiv -\infty$ or it is finite-valued 
on all of $I$ (and therefore convex).  The converse is immediate.\qed


\vfill\eject

\centerline{\bf  Appendix B. }
\smallskip
\centerline{\bf   Hessians of Plurisubharmonic Distributions. }
\medskip

 The decomposition
$$
\Symn\ =\ E\oplus S
\eqno{(B.1)}
$$
induces a decomposition of each $\Symn$-valued test function on $X$, and hence of 
each $\Symn$-valued distribution on $X$.
Applying this to $\Hess\,  u$, with $u\in \cd'(X)$, we have
$$
\Hess\,  u \ =\ (\Hess\,  u)_E + (\Hess\,  u)_S.
\eqno{(B.2)}
$$
\Lemma{B.1}  {\sl  If $u\in \PSH^{\rm dist}(X)$, then $(\Hess\,  u)_S$ is an $S$-valued measure
on $X$.}

\pf
Since the interior of $\plp$ in $S$ is non-empty, we may choose a basis $A_1,...,A_N\in 
{\Int \plp}$ for $S$.  the dual basis $A_1^*,...,A_N^*$ for $S$ will have the property that
$(\Hess\,  u)_S = u_1 A_1+\cdots+u_nA_N$ with each $u_j\in\cd'(X)$.
Given $\vf\in C^{\infty}_{\rm cpt}(X)$,
$$
u_j(\vf) \ =\  (\Hess\,  u)_S(\vf A_j).
\eqno{(B.3)}
$$
If $u$ is a \ppsh\ distribution, (B.3) implies that each $u_j\geq 0$ is a
non-negative measure. \qed\medskip

Note that using any basis for $\Symn$ (for example the standard  basis), 
$ (\Hess\,  u)_S$ will have measure coefficients.

\Lemma{B.2}  {\sl  Suppose $H$ is an $S$-valued measure on $X$.  Then there exists a measure $\|H\|\geq0$ and a   function  $\overrightarrow{H} : X\to S$
which is in $\lloc$ on $X$ with respect to the measure $\|H\|$, and $|\overrightarrow{H} (x)|=1$,
 $\|H\|$-a.e., such that
 $$
 H(\Phi)\ =\ \int_X \langle \overrightarrow{H} , \Phi\rangle \ \|H\|
 $$
for each $S$-valued test form $\Phi$ on $X$.  Also,  $\|H\|$ and $\overrightarrow{H}$
are unique.}
     
\pf
This is a standard fact about vector-valued measures.

\Theorem{B.3}  {\sl Suppose  $u\in\PSH^{\rm dist}(X)$ and abbreviate 
 $(\Hess\,  u)_S$  by $H_u$.  Then
 $$
  (\Hess\,  u)_S \ =\    \overrightarrow{H_u}  \ \|H_u\|
 \eqno{(B.4)}$$
   with $\|H_u\|\geq 0$ and $|\overrightarrow{H_u}(x)| =1$ for  $\|H_u\|$ a.a. $x\in X$.}


\vfill\eject

\centerline{\bf  Appendix C. }
\smallskip
\centerline{\bf   Convex Elliptic Sets in $\Symn$. }
\medskip

Suppose $F$ is an unbounded closed convex set in a finite dimensional inner product space
$(V, \langle \cdot,\cdot\rangle)$, and assume that $F$ has interior but $F\neq V$.
We can associate with $F$ two closed convex cones with vertex at the origin,
$\pp(F)$ and $\plp(F)$ which are polars of each other.

\medskip
\noindent
{\bf $\pp(F)$ --- The Ray Cone of $F$:}  Pick $a\in F$.  Consider the set of directions  $\oa u$
such that the ray from  $a$ in the direction $\oa u$ in contained in $F$. This coincides  with the compact subset
$$
\bigcap_{r>0} \smfrac 1 r \partial B_r \cap (F-a)
$$
of the unit sphere.  The cone on this compact set  is called the {\sl ray cone of $F$}
and is denoted by $\pp(F)$.  Since $F-a$ is convex, $\pp(F)$ is convex.
If $b\in F$ is any point in $F$, it is easy to see that the ray $\{a+tv:t\geq0\}$ is contained
in $F$ if and only if the ray $\{b+tv:t\geq0\}$ is contained in $F$.  That is,
$\pp(F)$ is independent of the choice of point $a\in F$.

\medskip
\noindent
{\bf $\plp(F)$ --- The  Cone of  Supporting Directions for $F$:}
For each non-zero $u\in V$ and each $\l\in \bbr$, consider the closed 
half-space
$$
H(u,\l)\ =\ \{v\in V: u\cdot v \geq \l\}
$$
If $F\ss H(u,\l)$ for some $\l\in \bbr$, then $u$ is a {\sl supporting direction vector for $F$}.
Let $\plp(F)$ denote the closure of the set of these supporting direction vectors.  
Obviously, $\plp(F)$ is a closed set of rays at the origin in $V$.  If $F\ss H(u,\l)$
and $F\ss H(u',\l')$ and $0\leq s\leq 1$, then it is easy to see that 
$
F\ss H(su+(1-s)u',s\l+(1-s)\l').
$
Thus $\plp(F)$ is convex.

\Prop{C.1}  {\sl  Suppose $F$ is an unbounded closed convex subset of $V$ with 
$\span F = V$ but $F\neq V$.  Then $\pp(F)$ and $\plp(F)$ are polars of each other
(with $\span \pp(F) =V$ and $\pp(F)\neq V$).}

\pf
Suppose $v\in\pp(F)$ and $u$ is a supporting direction vector.  Then for $a\in F$, the 
ray $\{a+tv : t\geq 0\} \ss F$ and there exists $\l\in\bbr$ with $F\ss H(u,\l)$.
Therefore,  $\l\leq \langle u, a+tv\rangle = \langle u, a\rangle +t \langle u, v\rangle$
for all $t\geq 0$ which implies that $\langle u, v\rangle\geq 0$.  This proves that each of
$\pp(F)$ and $\plp(F)$ is contained the the polar of the other.

Suppose $v$ is in the polar of $\plp(F)$, i.e.,  $\langle u, v\rangle \geq 0$ if $F\ss H(u,\l)$
for some $\l$.  Consider the ray  $\{a+tv : t\geq0\}$ through $a\in F$.
This ray is contained in $H(u,\l)$ since 
$\langle   a+tv, u\rangle = \langle a, u \rangle +t \langle v,u\rangle \geq \l+t\langle v,u\rangle\geq \l$
if $t\geq0$.  By the Hahn-Banach Theorem this ray must be contained in $F$.  Hence, $v\in \pp(F)$.
Thus $\pp(F)$ is the polar of $\plp(F)$.  The reverse follows from the bipolar theorem

\medskip
\noindent
{\bf The Edge of $F$.}
The set $E(F) = \{v\in V : \pm v \in \pp(F)\}$ consisting of those $v\in V$ such that the full affine line
$\{a+tv:t\in\bbr\}$ through $a\in F$ is contained in $F$, is called the {\sl linearity of $F$} 
or the {\sl edge of $F$}.  Set $S(F) \equiv E(F)^\perp$.
Then 
$$
F\ =\ E(F) \times (F\cap S(F))
$$
is a tube with base $F\cap S(F)$.  In this case the ray cone $\pp(F)$ is also a tube
$$
\pp(F) \ =\ E(F)\times (\pp(F)\cap S(F))
$$
with the same edge as $F$.

Note that $\span \pp(F) =V$ since $F$ is assume to have interior, but 
$$
\span \plp(F)\ =\ S(F).
$$

\medskip
\noindent
{\bf Convex Elliptic Sets.}
\Def{C.2}  A closed convex subset $F\ss\Symn$ which satisfies
\smallskip
\item{(1)}  \ \ $F+\cp\ \ss\ \ F$.
\smallskip
\item{(2)} \ \  $F$ can not be defined using the variables in a proper subspace $W\ss \rn$,
\smallskip
\noindent
will be called a {\sl convex elliptic set}.

\Prop{C.3}  {\sl A closed convex subset $F\ss\Symn$ is elliptic if and only if  its ray cone $\pp(F)$ is an elliptic cone.}

\pf It is easy to see that $F$ satisfies the positivity condition (1) if and only if $\pp(F)$ does.
It remains to show that $F$ can be defined using the variables in a proper subspace
$W\ss\rn$ if and only if the ray cone $\pp(F)$ can be defined using the variables in $W$

We must show that $\Sym(W)^\perp\ss F  \iff \Sym(W)^\perp \ss\pp(F)$.
One way is trivial.  For the other, suppose $\Sym(W)^\perp \ss\pp(F)$.
We may assume $0\in F$. Then 
$A\in \pp(F)$ if and only if the ray $\{tA : t\geq 0\}\ss F$.  Hence, $\Sym(W)^\perp\ss F$.\qed

\Cor{C.4}  {\sl   Suppose $F$ is a closed convex set in $\Symn$ with $F+\cp\ss F$.
Then $F$ cannot be defined  using fewer of the variables in $\rn$ if and only if each
$A\in\Int \plp(F)$ is positive definite.}


\vfill\eject

\centerline{\bf  Appendix D. }
\smallskip
\centerline{\bf   The Dirichlet Problem for Convex Elliptic Sets. }
\medskip

The main results of this paper carry over from elliptic cones   to convex elliptic sets $F$.
Suppose
$$
H\ =\ \{B\in \Symn : \langle A, B   \rangle \geq c\}
$$
is a supporting half-space for $F$ with $A\in \Int\plp(F)$.
By Corollary C.4, $A$ is positive definite.  Pick $B_0\in \partial H$, i.e., 
$ \langle A, B_0   \rangle  = c$.  Then the replacement for the mollifying condition
$\D_A u \geq 0$  is $\D_A u\geq \D_A B_0 =  \langle A, B_0   \rangle = c$.
The Mollifying Lemma \GG.2 remains valid for $F$-plurisubharmonic distributions.
The notion of  being u.s.c. $F$-\psh carries over in a straightforward manner.  
For both these concepts a function $u$ is of type $F$ if and only if it is of type $H$
for all supporting half-spaces $H$.
The key approximation property (7) is section 6 remains valid, with standard convolution
providing the proof.

The equivalent definitions of type $\wt F$ carry over from those of type ${\wt{\cp}}^+$.

Finally, the Dirichlet Problem is solvable in this context.
See Theorem 7.1 (Uniqueness) and Theorem 8.1 (Existence).
In the existence statement the boundary $\partial \Omega$ must be 
strictly $\pp(F)$-convex.

\Ex{\DD.1}  A simple but illuminating  example of a convex elliptic set is
$$
F\ \equiv\ \{ A\in \Symn :  A\geq0\ {\rm and\ } \det A \geq c\}
$$
for a constant $c\geq0$. One sees that $\pp(F) = \cp$.  The corresponding
classical equation is: $\det\{\Hess\, u\} = c$.

\Ex{\DD.2} A more interesting example is 
$$
F\ \equiv\ \{ A\in \Symn :  A\geq0\ {\rm and\ } {\rm Trace}\{ {\rm arctan} (A)\}  \geq  k\pi\}
\eqno{(\DD.1)}
$$
where $n=2k+1$ or $2k+2$.  The corresponding equation
$$
{\rm Im} \{ \det(I+iA)\}\ =\ 0
\eqno{(\DD.2)}
$$
 arises in the study of Special Lagrangian submanifolds, and the Dirichlet problem
 for (\DD.2) was studied in depth by Caffarelli, Nirenberg and Spruck [CNS].
In fact the locus of (\DD.2) has $k$ connected components and [CNS] treats 
only the ``outermost''  component, which corresponds to the boundary of the
set $F$ defined in (\DD.1). In [CNS] the authors show that 
$$
\pp(F) \ =\ 
\cases{
\cp\qquad {\rm if \ } n {\rm \ is\ odd}  \cr
{\cal Q}\qquad {\rm if \ } n {\rm \ is\ even}  \cr
}
$$
where $\Int {\cal Q}$ is the  component of the set 
$\{A\in\Symn: \sigma_{n-1}(A) >0\}$ which contains the identity $I$
(and $\sigma_{n-1}$ denotes the $(n-1)$st elementary symmetric function).
In [HL$_4$] existence and uniqueness of continuous solutions to the
 Dirichlet Problem are established for all branches of the equation (\DD.2).
However, the smoothness of the solutions for smooth boundary data remains
largely open (see [Y] however).


\vfill\eject

\centerline{\bf  Appendix \WW. }
\smallskip
\centerline{\bf   Elliptic $\MA$-operators / G\aa rding-Hyperbolic Polynomials on Sym$^2(\rn)$. }
\medskip

\def\MA{\rm MA}

For each polynomial $P$ on the vector space $\Symn$ consider the associated (non-linear) 
partial differential operator defined by ${\bf P}(f) = P(\Hess f)$.  If $P$ is the determinant,
the associated operator is the real Monge-Amp\`ere operator.

\Def{\WW.1}  Let $M$ be a homogeneous polynomial of degree $m$ on $\Symn$.
Suppose that the identity is a hyperbolic direction for $M$ in the sense of G\aa rding [G]. 
That is, suppose that for each $A\in\Symn$, the polynomial $p_A(t)=M(tI+ A)$ has exactly $m$ real zeros on $\bbr$, and that $M(I)=1$.  Then the operator 
$$
\M(f)\ = \ M(\Hess f)
\eqno{(\WW.1)}
$$
will be called an {\sl $\MA$-operator}, and the polynomial $M$ will be called an {\sl
$\MA$-polynomial.}
\medskip

G\aa rding's beautiful theory of hyperbolic polynomials states that the set
$$
\G(M)\ =\ \{A\in\Symn: M(tI+A)\neq 0\ \ {\rm for\ } t\geq0\}
\eqno{(\WW.2)}
$$
is an open convex cone in $\Symn$ equal to the connected component of $\{M>0\}$ containing $I$.
The closed convex cone
$$
\pp(M)\ =\ \{A\in\Symn: M(tI+A)\neq 0\ \ {\rm for\ } t>0\}
\eqno{(\WW.3)}
$$
is the closure of $\G(M)$.  Moreover, 
$$
\partial \pp(M)\ =\ \{A\in\Symn:  M(A) =0 {\ \rm but\ } M(tI+A)\neq 0\ \ {\rm for\ } t>0\}.
$$
Let $\plp$ denote the polar cone to $\pp(M)$.
\smallskip
\noindent
The Positivity Condition  on $\pp$ (from \S 3) can be stated in several equivalent ways in terms of $M$:
\smallskip

1)\ \ \ $M(tI+A)\neq 0$ for all $t>0$ and $A>0$ \ (i.e., $\Int \cp \ss \pp(M)$).

\smallskip

1)$'$\ \ \ $M(tI+P_e)\neq 0$ for all $t>0$ and all unit vectors $e$  \ (i.e., $P_e \in\pp(M)$ for all $e$).
\smallskip

\noindent
The Completeness Condition  on $\pp$ can also be stated in several equivalent ways in terms of $M$:
\smallskip

2)\ \ \ $M(tI-P_e)$  has a strictly positive zero  for each unit vector $e$ \ (i.e., $P_e \notin$
  the edge 
  
   \qquad of  $(\pp(M))$ for each $e$).

\smallskip

2)$'$\ \ \  For each $A$ and each unit vector $e\in \rn$,  $M(tP_e+A)$  is non-constant in $t$.
\smallskip
  \Prop{\WW.2}  {\sl  The cone $\pp(M)$ defined by an $\MA$-polynomial is elliptic if and only if for each unit vector $e\in \rn$,\smallskip
  
  a)\ \ \ $M(I+sP_e)$ is not $\equiv 1$, and \smallskip
  
  b)\ \ \  $M(I+sP_e)>0$  for $s>0$.}

\medskip

The linearization of the non-linear operator $\M$ at a point $x$ and a function $f$ is
$$
L(g)\ =\ L_A(g) \ =\ \smfrac{d}{dt} M(A+tH)\bigr|_{t=0}\ =\ m\overline M (H,A,...,A)
\eqno{(\WW.4)}
$$
where $A=\Hess_x f$, $H=\Hess_x g$, and $\overline M$ is the completely polarized form of $M$.
The linear  functional $L_A$ on $\Symn$ determines a unique element $\wt A\in \Symn$
such that $L_A(H)=\langle H, \wt A  \rangle$, and $L_A$ is elliptic if and only if $\wt A$ is positive definite.  (If $M$ is homogeneous of degree $m$, then $A\mapsto \wt A$ is homogeneous of
degree $m-1$.)

The next result helps to justify the terminology ``elliptic cone'' introduced in \S 3.

\Theorem {\WW.3}  {\sl  Suppose $\M$ is an $\MA$-operator.  Then $\M$ is elliptic at each $f$, $x$ with 
$f$ strictly \ppsh \ (equivalently at each $\Hess_x f = A\in \pp(M)$) if and only if the cone $\pp(M)$ is an elliptic cone.
}

\pf
Suppose $\M$ is an elliptic operator at each $A\in {\rm Int}\pp(M)$, i.e., the symmetric form
$\wt A \in \Sym(\rn)$ defined by 
$$
{d\over dt} M(A+tH)\biggr|_{t=0}\ =\ \langle \wt A, H   \rangle
$$
is positive definite.  By G\aa rding's inequality [G] 
$$
\langle \wt A, H   \rangle \ =\ m \overline M(H; A,...,A)\ >\ 0
$$
if $H\in {\rm Int}\pp(M)$.  Hence $\wt A\in\plp(M)$ and $\wt A$ is positive definite.

Since $\wt A$ is positive definite, for each $e\in\rn$ with  $|e|=1$, we have that 
$0< \langle \wt A, P_e   \rangle = {d\over dt} M(A+tP_e)\bigr|_{t=0}$.  By the same argument,
$$
 {d\over dt} M(A+tP_e)\ >\ 0\qquad {\rm if}\ \ A+tP_e\in {\rm Int}\pp(M).
$$
This implies that $M(A+tP_e)>0$ for all $t>0$, i.e., the ray
$\{A+tP_e) : t\geq 0\} \subset  {\rm Int}\pp(M)$ for all $A\in  {\rm Int}\pp(M)$.  Equivalently,
$P_e\in \pp(M)$.  This proves the Positivity Condition for $\pp(M)$.

Suppose  now that $\pp(M)$  is elliptic.
Then $P_e\notin E$, the edge of $\pp(M)$.
By Theorem 3, p. 962 in [G], $M(A+tP_e)$ is not constant in $t$.
Suppose $A\in {\rm Int}\pp(M)$. By   the Positivity Condition there exists  $\l>0$ such that 
$A+tP_e \in {\rm Int}\pp(M)$ for all $t\in (-\l,+\infty)$.
Define $g(t) =M(A+tP_e)^{1\over m}$ on $(-\l,+\infty)$.
Then $g>0$, and by [G], $g$ is concave on $(-\l,+\infty)$.  As noted $g$
is not constant.  A concave function on $(-\l,+\infty)$ which is $>0$  and non-constant,
such as $g$,  must be strictly increasing.  Therefore, $g(t)^m$ is also strictly increasing.
This proves that $\langle \wt A, P_e   \rangle = {d\over dt} M(A+tP_e)\bigr|_{t=0} >0$ for each $P_e$. 
Therefore, $\wt A$ is positive definite.
\qed

\vfill\eject
\noindent
{\bf Examples:}  The basic examples are given by the determinant.  There are four cases
corresponding to $\bbr, \bbc, \bbh$ and $\bbo$.  \medskip

\noindent
1.\ \  The determinant on $\Symn$.     \medskip

\noindent
2.\ \  The  determinant on \ ${\rm Herm}_\bbc\Sym(\bbc^n)  \ \ss\ \Sym(\bbr^{2n})$.     \medskip

\noindent
3.\ \  The  determinant on \ ${\rm Herm}_\bbh\Sym(\bbh^n)  \ \ss\ \Sym(\bbr^{4n})$.     \medskip

\noindent
4.\ \  The  determinant on \ ${\rm Herm}_\bbo\Sym(\bbo^2)  \ \ss\ \Sym(\bbr^{16})$.     \medskip

\medskip

The quaternionic case is  perhaps best understood as a polar action [DK].  Namely,
Sp$_n$ acts on ${\rm Herm}_\bbh\Sym(\bbh^n)$ with cross-section given by the space
$D$ of diagonal matrices.  The polynomial $\l_1\cdots \l_n$ on $D$ extends to an 
Sp$_n$-invariant polynomial, det, on ${\rm Herm}_\bbh\Sym(\bbh^n)$ (cf. [AV]).

In each of these cases the inhomogeneous equation has been been treated:  the real case
by Taylor-Rauch, the complex case by Bedford-Taylor, the quaternionic case by Alesker-Verbitsky
and the octonian case also by Alesker-Verbitsky.

Certain versions of the inhomogeneous Monge-Amp\`ere can be treated by the methods
in [HL$_{4,5,6}$]. For example one can insert a function  $f(x,u)$ with $f_u\geq0$.
One can also address all other branches of the determinant in this inhomogeneous form.
\medskip

Finally, if a polynomial $M$ as above is hyperbolic in the direction $I\in \Sym(\bbr^N)$,
then $M^{(k)}(A) \equiv \overline M(I,...,I;A,...,A)$ with $A$ inserted into $k$ slots, is also hyperbolic
in the direction $I$.  Thus the elementary symmetric functions provide additional examples
in all of these cases.

 \vfill\eject



\centerline{\bf References}

\vskip .2in

\noindent
\item{[Al]}   S. Alesker,  {\sl  Non-commutative linear algebra and  plurisubharmonic functions  of quaternionic variables}, Bull.  Sci.  Math., {\bf 127} (2003), 1-35. also ArXiv:math.CV/0104209.  

\smallskip
 
\noindent
\item{[AV]}    S. Alesker and M. Verbitsky,  {\sl  Plurisubharmonic functions  on hypercomplex manifolds and HKT-geometry}, arXiv: math.CV/0510140  Oct.2005

\smallskip

\noindent
\item{[BT]}   E. Bedford and B. A. Taylor,  {The Dirichlet problem for a complex Monge-Amp\`ere equation}, 
Inventiones Math.{\bf 37} (1976), no.1, 1-44.

\smallskip

 \item{[B]}  H. J. Bremermann,
    {\sl  On a generalized Dirichlet problem for plurisubharmonic functions and pseudo-convex domains},
          Trans. A. M. S.  {\bf 91}  (1959), 246-276.
\smallskip

\noindent
\item{[C]}   M. G. Crandall, {\sl Viscosity solutions: a primer}, pp. 1-37 in Viscosity Solutions and Applications, Eds.  I. C. Dolcetta and P. L. Lions, 
Springer Lecture Notes in Math. 1660, 1995.

 \smallskip

\noindent
\item{[CIL]}   M. G. Crandall, H. Ishii and P. L. Lions {\sl
User's guide to viscosity solutions of second order partial differential equations},  
Bull. Amer. Math. Soc. (N. S.) {\bf 27} (1992), 1-67.

 \smallskip

\noindent
 \item{[CNS]}   L. Caffarelli, L. Nirenberg and J. Spruck,  {\sl
The Dirichlet problem for nonlinear second order elliptic equations, III: 
Functions of the eigenvalues of the Hessian},  Acta Math.
  {\bf 155} (1985),   261-301.

 \smallskip

\noindent
\item{[DK]}   J. Dadok and V. Katz,   {\sl Polar representations}, 
J. Algebra {\bf 92} (1985) no. 2,
504-524.

\smallskip

\noindent
\item{[D]}  J.-P. Demailly, { Complex Analytic and Differential Geometry},
  e-book at Institut Fourier, UMR  5582 du CNRS,
   Universit\'e de Grenoble I, Saint-Martin d'H\`eres, France:
   can be found at http://www-fourier.ujfgrenoble.fr/~demailly/books.html.
 \smallskip

\noindent
\item{[G]}   L. G\aa rding, {\sl  An inequality for hyperbolic polynomials},
 J.  Math.  Mech. {\bf 8}   no. 2 (1959),   957-965.

 \smallskip

   \noindent 
\item {[HL$_1$]} F. R. Harvey and H. B. Lawson, Jr, {\sl Calibrated geometries},  Acta Mathematica 
{\bf 148} (1982), 47-157.

 \smallskip

   \noindent 
 \item{[HL$_2$]} F. R. Harvey and H. B. Lawson, Jr, {\sl An introduction to potential theory in calibrated geometry}, Stony Brook Prerprint, 2006.  ArXiv:0710.3920

 \smallskip

   \noindent 
\item {[HL$_3$]} F. R. Harvey and H. B. Lawson, Jr, {\sl Duality of positive currents and plurisubharmonic functions in calibrated geometry},  Amer. J. Math. (to appear).  ArXiv:math.0710.3921.
 \smallskip

   \noindent 
\item {[HL$_4$]} F. R. Harvey and H. B. Lawson, Jr, {\sl  Dirichlet Duality and the Nonlinear Dirichlet Problem}, Stony Brook Prerprint, 2008. ArXiv:0710.3991

 \smallskip

   \noindent 
\item {[HL$_5$]} F. R. Harvey and H. B. Lawson, Jr, {\sl Potential Theory for Second Order Partial Differential Equations on $\rn$},  

 \smallskip

   \noindent 
\item {[HL$_6$]} F. R. Harvey and H. B. Lawson, Jr, {\sl  Potential Theory for Second Order Partial Differential Equations on Manifolds},  

 \smallskip

\noindent
\item{[HW$_1$]} F. R. Harvey,  R. O. Wells, Jr.,  {\sl Holomorphic approximation and hyperfunction 
theory on a $C^1$ totally real submanifold of a complex manifold},
  Math.  Ann. {\bf 197} (1972),  287-318.

 \smallskip

\noindent
\item{[HW$_2$]}  F. R. Harvey,  R. O. Wells, Jr.,  {\sl Zero sets of non-negatively strictly plurisubharmonic
functions},
  Math.  Ann. {\bf 201} (1973),  165-170.

 \smallskip

   \noindent
\item{[J]}    R. Jensen,    {\sl  Uniqueness criteria for viscosity solutions of fully nonlinear 
elliptic partial differential  equations},    Indiana Univ. Math. J. {\bf 38}  (1989),   629-667.

\smallskip

 \noindent 
\item {[M]} J. Milnor, {Morse Theory}, Annals of Math. Studies no. {\bf 51}, Princeton University Press,
 Princeton, N.J.,  1963.
 \smallskip

\smallskip

\item {[RT]} J. B. Rauch and B. A. Taylor, {\sl  The Dirichlet problem for the 
multidimensional Monge-Amp\`ere equation},
Rocky Mountain J. Math. {\bf 7}    (1977), 345-364.

\smallskip

\item {[R]} R. Richberg, {\sl  Stetige streng pseudokonvexe Funktionen},
  Math. Ann. {\bf 175}    (1968), 257-286.

\smallskip

\item {[S]}  Z. Slodkowski, {\sl  The Bremermann-Dirichlet problem for $q$-plurisubharmonic functions},
Ann. Scuola Norm. Sup. Pisa Cl. Sci. (4)  {\bf 11}    (1984),  303-326.

\smallskip

\item {[W]}   J.  B. Walsh,  {\sl Continuity of envelopes of plurisubharmonic functions},
 J. Math. Mech. 
{\bf 18}  (1968-69),   143-148.

\smallskip

\item {[Y]}  Yu Yuan,  {\sl A priori estimates for solutions of fully nonlinear special lagrangian equations},
 Ann Inst. Henri Poincar\'e  
{\bf 18}  (2001),   261-270.

\end